\documentclass[12pt]{article}
\usepackage{amsfonts}
\usepackage{amssymb}
\usepackage{amscd}
\usepackage{amsmath}
\usepackage{epsfig}
\usepackage{graphicx, subfigure}
\usepackage{latexsym}
\usepackage{pst-3dplot}
\newlength{\dinwidth}
\newlength{\dinmargin}
\newtheorem{theorem}{Theorem}

\newtheorem{proposition}{Proposition}
\newtheorem{corollary}{Corollary}
\newtheorem{remark}{Remark}
\newtheorem{lemma}{Lemma}
\newtheorem{example}{Example}

\def \surf{{\cal L}}
\def\rank{\mathrm{rank}}

\begin{document}

\title{Algebro-geometric solutions of the Schlesinger systems and the Poncelet-type polygons in higher dimensions}
\author{Vladimir Dragovi\'c$^1$ and Vasilisa Shramchenko$^2$}

\date{}

\maketitle

\footnotetext[1]{Department of Mathematical Sciences, University
of Texas at Dallas, 800 West Campbell Road, Richardson TX 75080,
USA. Mathematical Institute SANU, Kneza Mihaila 36, 11000
Belgrade, Serbia.  E-mail: {\tt
Vladimir.Dragovic@utdallas.edu}--the corresponding author}

\footnotetext[2]{Department of mathematics, University of
Sherbrooke, 2500, boul. de l'Universit\'e,  J1K 2R1 Sherbrooke, Quebec, Canada. E-mail: {\tt Vasilisa.Shramchenko@Usherbrooke.ca}}

\begin{abstract}
A new method to construct algebro-geometric solutions of rank two
Schlesinger systems is presented. For an elliptic curve represented as a ramified double covering of $\mathbb{CP}^1$,
a meromorphic differential is constructed with the following property: the common projection of its two  zeros on the base of the covering, regarded as a
 function of the only moving branch point of the covering, is a solution of a Painlev\'e VI equation. This differential provides an invariant
formulation of a classical Okamoto transformation for the Painlev\'e VI equations. A generalization of this differential to hyperelliptic curves is also constructed. In this case, positions of zeros of the differential provide part of a solution of the multidimensional Garnier system. The corresponding solutions of the rank two  Schlesinger systems associated with elliptic and hyperelliptic curves are constructed in terms of this differential. The initial data for construction of the meromorphic differential include a point in the
Jacobian of the curve, under the assumption that this point has nonvariable
coordinates with respect to the lattice of the Jacobian while the
branch points vary. It appears that the cases where the
coordinates of the point are rational correspond to periodic
trajectories of the billiard ordered games associated with $g$ confocal quadrics in $g+1$
dimensional space. This is a generalization of a situation studied by
Hitchin, who related algebraic solutions of a Painlev\'e VI
equation with the Poncelet polygons.

\vskip 1cm



\end{abstract}

\newpage

\tableofcontents

\bigskip
\section{Introduction}

Two approaches to algebro-geometric solutions of the rank two Schlesinger
systems, both based on the Riemann-Hilbert problem technique were suggested in
\cite{Deift} and \cite{KiKo}. Both approaches resulted in a $2g$-parametric family of  explicit solutions written in terms of theta-functions associated with hyperelliptic curves of genus $g$.
We are presenting an essentially different method of constructing algebro-geometric
solutions of the Schlesinger systems. Our  point of departure is an observation by Hitchin
that the Poncelet polygons in the plane are closely related to algebraic
solutions of a certain classical Painlev\'e VI equation, see
\cite{Hitchin}.

The Painlev\'e VI equation is the second order ordinary differential equation
\begin{equation}
\label{Painleve}
\frac{d^2 y}{dx^2} = \frac{1}{2} \left(  \frac{1}{y} + \frac{1}{y-1}+ \frac{1}{y-x}  \right) \left( \frac{dy}{dx} \right)^2 - \left( \frac{1}{x} + \frac{1}{x-1} + \frac{1}{y-x} \right) \frac{dy}{dx}
\end{equation}
\begin{equation*}
 + \frac{y(y-1)(y-x)}{x^2(x-1)^2}\left( \alpha +\beta\frac{x}{y^2} + \gamma\frac{x-1}{(y-1)^2} + \delta\frac{x(x-1)}{(y-x)^2}  \right)
\end{equation*}
with parameters $\alpha, \beta, \gamma, \delta \in \mathbb C$.  Historically this equation was derived by Gambier \cite{Gambier} and appeared in 1907 \cite{Fuchs} in R. Fuchs'
study of isomonodromic deformations
of the Fuchsian linear system
\begin{equation*}
\frac{d\Phi}{du} = A(u) \Phi, \qquad u \in {\mathbb C},
\end{equation*}
for a $2\times 2$ matrix function $\Phi(u)$ defined in the complex
plane, where the matrix $A\in sl(2,{\mathbb C})$ is of the form:
\begin{equation*}
\label{A} A(u) = \frac{A^{(1)}}{u} + \frac{A^{(2)}}{u-1} +
\frac{A^{(3)}}{u-x}.
\end{equation*}
The  residue-matrices $A^{(i)}\in sl(2,{\mathbb C})$, $i=1,2,3$, are
$u$-independent and, as a consequence  of the isomonodromy condition, satisfy the Schlesinger system:
\begin{equation}
\label{Schlesinger}
\frac{dA^{(1)}}{dx} = \frac{[A^{(3)}, A^{(1)}]}{x}; \qquad
\frac{dA^{(2)}}{dx} = \frac{[A^{(3)}, A^{(2)}]}{x-1}; \qquad
\frac{dA^{(3)}}{dx} = -\frac{[A^{(3)}, A^{(1)}]}{x} -
\frac{[A^{(3)}, A^{(2)}]}{x-1}.
\end{equation}
The last equation is equivalent to  $A^{(\infty)}:= -A^{(1)} -A^{(2)}-A^{(3)} = const.$

Using the freedom of  global conjugation of the residue matrices by a constant invertible matrix, one may assume  $A^{(\infty)}$ to be diagonal.  Then the entry $A_{12}(u)$ of the matrix $A(u)$ has only one zero in the complex plane and  is of the form:
\begin{equation}
\label{A12}
A_{12}(u) = \kappa\frac{(u-y)}{u(u-1)(u-x)}
\end{equation}
with some constant $\kappa\in{\mathbb C}.$ The position of the zero of $A_{12}(u)$ as a function of the position $x$ of the pole is the function $y(x)$ which satisfies
the Painlev\'e VI equation (\ref{Painleve}) with parameters
\begin{equation*}
\alpha = \frac{(t_1-1)^2}{2}, \qquad \beta=-\frac{t_2^2}{2}, \qquad \gamma=\frac{t_3^2}{2}, \qquad \delta = \frac{1}{2} - \frac{t_4^2}{2}.
\end{equation*}
Here, $\frac{t_i}{2}$ and $-\frac{t_i}{2}$ are the eigenvalues of the residue-matrix $A^{(i)}$; these quantities are integrals of motion of system (\ref{Schlesinger}).

In this paper we consider Painlev\'e VI equation (\ref{Painleve})  with the following parameters, which corresponds to a solution of the Schlesinger system with eigenvalues $\frac{t_i}{2} = \pm\frac{1}{4}$:
\begin{equation}
\label{constants}
\alpha=\frac{1}{8}, \qquad \beta=-\frac{1}{8}, \qquad \gamma=\frac{1}{8}, \qquad \delta=\frac{3}{8}.
\end{equation}
This equation is particular in the family (\ref{Painleve}) of
Painlev\'e VI equations since its general solution can be written
in terms of known functions. Note that other Painlev\'e VI equations with this property are described in \cite{Okamoto} and \cite{Manin}; and in  \cite{LT} and \cite{Iwasaki} all Painlev\'e VI admitting algebraic solutions are classified.
In \cite{Dubrovin} the Painlev\'e VI equations (\ref{Painleve}) appeared
in the context of Frobenius manifolds, and in Appendix E
algebraic solutions were constructed for a class of equations
(the class IV according to \cite {LT}, p. 44, which contains the case of
parameters (\ref{constants})).
%
The  general
solution of the Painlev\'e VI equation \eqref{Painleve}, \eqref{constants}  in terms of Jacobi theta-functions was given in
\cite{Hitchin_Painleve}. Later on, the same solution was rewritten
in a different, more compact form in \cite{BK}. In this paper, we
find yet another interpretation of the general solution to the
Painlev\'e VI equation with coefficients (\ref{constants}) in
terms of position of zeros of a certain Abelian differential of
the third kind on the associated elliptic curve. Our approach is
motivated by the work of Hitchin \cite{Hitchin} which we now
briefly outline.

For two plane conics, a Poncelet polygon is a closed trajectory
inscribed into one of the conics, called the boundary, and
circumscribed about the other one, called the caustic. Two given
conics generate a one-parameter family of conics inscribed in the
four common tangents of the given conics, a tangential pencil of
conics. Suppose now that a tangential pencil of conics is given
and one of the conics of the pencil is fixed, the caustic. One can
associate an elliptic curve with this data (pencil of conics, a
caustic from the pencil), see for example \cite{GriffithsHarris},
\cite{DR2011}. One can ask for which values of the parameter of
the pencil, the pair caustic-boundary, where the caustic has been
fixed above, and the boundary corresponds to the value of the
pencil parameter, is such that there exists a Poncelet polygon for
it. For the given caustic and boundary, the Poncelet theorem
states that if there exists one Poncelet polygon with $k$
vertices, then there are infinitely many such polygons with the
same number $k$ of vertices and any point of the boundary is a
vertex of such a polygon (see Fig. \ref{fig:sedmouglovi} for
$k=7$). For a given pencil of conics with a fixed caustic, the
Poncelet polygons appear when the boundary conic corresponds to a
point of a finite order on the associated elliptic curve.

Let the elliptic curve $\mathcal L$ be given by the following
equation:
$$
\mathcal L:\;\; v^2=u(u-1)(u-x).
$$
Denoting by $\mathcal A$ the corresponding Abel map based at the ramification point at infinity and by
$(1\;\; \mu)$ the vectors generating the lattice of the Jacobian of the elliptic curve, a point $Q_0 \in \mathcal L$ is of the
order $k\in{\mathbb N}\;$ if
$\; k\mathcal A(Q_0)=  n+m\mu,$ {with some
integer constants} $ m$ and $n$. In other words, $Q_0$ is of the finite order if
\begin{equation}
\label{eq:rational}
 \mathcal A(Q_0)=  c_1+c_2\mu,
\end{equation}
with some rational constants $ c_1, c_2 \in{\mathbb Q}$.

However, if the caustic is not fixed but varies within the pencil,
then the associated elliptic curve also varies by changing the
position of the branch point $x$.

Hitchin showed in \cite{Hitchin} that for any point of a finite
order $k$ different from a ramification point on the elliptic curve $\mathcal L$, and thus for any Poncelet polygon of length $k\geq 3$, one can construct an algebraic solution of the Painlev\'e VI equation with parameters (\ref{constants}). Following Hitchin, for such a point $Q_0$ of finite order $k$ on the elliptic curve $\mathcal L$, one can construct a differential $\Omega$
on $\mathcal L$ with simple poles at $Q_0, Q_0^*$, where
$*$ is the elliptic involution on $\mathcal L$, with
the property that the common $u$-coordinate of the two zeros of the differential considered as
a function of the branch point $x$ is an algebraic solution of the
Painlev\'e VI equation \eqref{Painleve}, \eqref{constants}.

More precisely, as our first remark and an initial step in our construction, we
describe this differential explicitly as follows. Suppose some
canonical basis $a, b$ of homology is chosen on $\mathcal L$ and
denote by $\omega$ the holomorphic differential on the curve
normalized by the condition $\oint_a\omega=1.$ Let $\Omega_{Q_0,
Q_0^*}$ be the differential of the third kind with poles at $Q_0$
and $Q_0^*$ with residues $+1$ and $-1$, respectively, normalized
to have vanishing $a$-period. Then, we have the
differential $\Omega$ in the form
\begin{equation}
\label{Omega_intro}
\Omega = \Omega_{Q_0,Q_0^*} - 4 \pi {\rm i} c_2 \omega.
\end{equation}
This new form of the differential $\Omega$ turns out to be the cornerstone of further generalizations. Alternatively, we define
$\Omega$ as an Abelian differential of the third kind with the
poles as above and the $a$-period equal to $- 4 \pi {\rm i} c_2$.
Then, due to the position of the poles, the $b$-period of $\Omega$
is equal to $4 \pi {\rm i} c_1.$ This differential has two zeros
on the curve $\mathcal L$. The zeros are related by the elliptic
involution. Now $x$ is considered an independent variable and the point $Q_0$ on the varying elliptic curve is defined by condition (\ref{eq:rational}) with constant values of $c_1$ and $c_2$.  Then the common $u$-coordinate of two zeros of $\Omega$ as a function of $x$ is an
algebraic solution to the Painlev\'e VI equation with parameters
(\ref{constants}) that is shown in \cite{Hitchin} to correspond to
a Poncelet polygon of length $k$.

Our second remark is that the condition on rationality of $c_1,
c_2$ can be relaxed. If we assume that in the condition of type
(\ref{eq:rational})  $c_1$ and  $c_2$ are  arbitrary complex
constants not being simultaneously equal to half-integers, then
the zeros of the differential  $\Omega$ are still solutions of the
Painlev\'e VI equation with parameters (\ref{constants}).
Moreover, the differential $\Omega$ provides an invariant
interpretation of a classical Okamoto transformation,  as long as
the point $Q_0$ satisfies condition (\ref{eq:rational}) with the
constants $c_1$ and $c_2$
independent of variations of the branch point $x$ and such that $(c_1, c_2)$ is not a half-integer vector.

The described variation of the pair, the elliptic curve $\mathcal L$ and a point $Q_0$ on it, in Hitchin's case \cite{Hitchin} can be interpreted as a simultaneous
deformation of the boundary and the caustic such that the rotation
number (see for example \cite{KT1991}, \cite{Tab2005book},  and
\cite{DR2011} the definition of the rotation number) remains fixed
and rational.
In our case, this variation corresponds to a simultaneous deformation of the
boundary and the caustic such that the rotation number remains
constant, but not necessarily rational.

We then use the differential $\Omega$ to construct a solution of the Schlesinger system (\ref{Schlesinger}), starting from the fact that both, our differential and the term $A_{12}(u)$, see (\ref{A12}), of the connection matrix $A(u)$, have their only zero at the solution to the Painlev\'e VI equation with parameters (\ref{constants}).

Our next observation is that the described interrelation between
the Poncelet polygons, the Painlev\'e VI equation, and elliptic
curves admits a natural generalization to the hyperelliptic
curves. The role of the Painlev\'e VI equation is then played by
the multidimensional Garnier system described in \cite{Jap},
p. 169. Instead of the  Poncelet polygons inscribed in a conic, we
then obtain billiard trajectories bouncing off a family of
quadrics in a multidimensional space.

More precisely, consider a more general situation of the matrix linear system for a $2\times 2$ matrix function $\Phi(u)$ defined in the Riemann sphere
\begin{equation}
\label{linsys_intro}
\frac{d\Phi}{du} = A(u) \Phi, \qquad u \in {\mathbb CP}^1
\end{equation}
 with the matrix $A\in sl(2,{\mathbb C})$ having $2g+1$ simple poles at $u_1, u_2, \dots,
u_{2g+1}$ and at the point at infinity:
\begin{equation}
\label{Ah}
A(u) = \sum_{j=1}^{2g+1} \frac{A^{(j)}}{u-u_j}.
\end{equation}
The corresponding isomonodromic deformation equations form the
Schlesinger system of partial differential equations for the residue-matrices $A^{(j)}\in
sl(2,{\mathbb C})$ with respect to positions $u_i$ of poles of $A$
as independent variables:

\begin{equation}
\label{hSchlesinger_intro}
\frac{\partial A^{(j)}}{\partial u_k} =
\frac{[A^{(k)}, A^{(j)}]}{u_k-u_j}; \qquad\qquad  \frac{\partial
A^{(k)}}{\partial u_k} = -\sum_{j\neq k}\frac{[A^{(k)},
A^{(j)}]}{u_k-u_j},
\end{equation}
where  $ A^{(\infty)}:=-A^{(1)} -\dots-A^{(2g+1)}  = const.$

This Schlesinger system is naturally associated with the family of hyperelliptic curves defined by the equation
\begin{equation}
\label{hyper_intro} v^2=(u-u_1)\cdots(u-u_{2g+1}),
\end{equation}
where the varying branch points $u_1, \dots, u_{2g+1}$ are given by the positions of the poles in (\ref{Ah}).

Generalizing the pattern we have developed in the elliptic case to
the hyperelliptic curves, we  construct a meromorphic
differential $\Omega$, the hyperelliptic analogue of (\ref{Omega_intro}). Namely, assume a canonical homology basis $\{a_1, b_1, \dots, a_g, b_g\}$ is chosen on a hyperelliptic curve from the family (\ref{hyper_intro}),
 let  $c_1,c_2 \in {\mathbb C}^g$ be two constant vectors
 and consider a point $z_0$ in the Jacobian of the curve such that
$$
z_0: = c_1 + {\mathbb B}c_2.
$$
Here $\mathbb B$ is the Riemann matrix of the curve and the matrix $({\rm I}\,\,\, \mathbb B)$
generates the lattice $\Lambda$ of the Jacobian of the curve. The
Jacobi inversion of $z_0$ is a positive  divisor  $Q_1+\dots +
Q_g$ of degree $g$ on the hyperelliptic curve. Let us assume that the
points $Q_j$ are all distinct and none of them coincides with a
branch point of our hyperelliptic curve and that no two points of
the divisor are paired by the hyperelliptic involution
interchanging the sheets of the covering of the $u$-sphere. As the
main tool of our construction, we define the following
differential:
\begin{equation*}
\Omega(P ) = \sum_{j=1}^g \Omega_{Q_j Q_j^*}(P ) - 4\pi {\rm i}\, c_2^t{\bf \omega}(P ).
\end{equation*}
Here $Q_j^*$ is the hyperelliptic involution of the point $Q_j$ and $\Omega_{Q_j Q_j^*}$ is the Abelian differential of the third kind with simple poles at $Q_j$ and $Q_j^*$ of residue $+1$ and $-1$, respectively, normalized by vanishing of all its $a$-periods; and $\omega$ is a vector of holomorphic differentials $\omega_j$ normalized by $\oint_{a_j}\omega_i=\delta_{ij}$.  In other words, $\Omega$ is a meromorphic differential with the given pole structure such that its $a_k$-period is $-4\pi{\rm i}c_{2k}$ while its $b_k$-period is $4\pi{\rm i}c_{1k}$ with $c_{1k}$ and $c_{2k}$ standing for the $k$th components of the vectors  of constants $c_1$ and $c_2.$ This differential has $4g-2$ zeros paired by the hyperelliptic involution.

Having constructed the differential $\Omega$ on one hyperelliptic curve from the family (\ref{hyper_intro}), we start varying the curve by allowing small variations of positions of the branch points $\{u_j\}_{j=1}^{2g+1}$. Note that the lattice $\Lambda$ will vary accordingly, however we define the point $z_0$ by the requirement that the vectors $c_1, c_2$, its coordinates with respect to the varying lattice $\Lambda$, stay constant.

We then construct a $2g$-parametric family, having the vectors  $c_1$ and $c_2$ as parameters, of solutions to the Schlesinger system  (\ref{hSchlesinger_intro})  in terms of the  differential $\Omega$ by finding a hyperelliptic analogue of our solution to the four-point Schlesinger system (\ref{Schlesinger}) corresponding to the family of elliptic curves.

As a consequence of our construction and Corollary $6.2.2$ in \cite{Jap}, the $2g-1$ functions of the branch points $\{u_j\}_{j=1}^{2g+1}$ given by the $u$-coordinates of the pairs of zeros of $\Omega$ give solutions to the multidimensional Garnier system of the form given in \cite{Jap}.

We note that the starting point for our construction is a point $z_0$ in the Jacobian of a hyperelliptic curve with coordinates
 $c_1,c_2 \in {\mathbb C}^g$ with respect to the lattice $\Lambda$ generated by $({\rm I}\,\,\, \mathbb B)$.
In the case when these coordinates are rational, $c_1,c_2 \in {\mathbb Q}^g$, which can be restated as follows
$$
z_0\in \frac{1}{k}\Lambda, \qquad k\in \mathbb N, \;\;\;k\geq3,
$$
one can associate such a  point $z_0$  with the  periodic trajectories of
the so-called billiard ordered game for $g$ confocal quadrics in the $(g+1)$-dimensional space. The periodic trajectories of the billiard ordered
game are analogues in the $(g+1)$-dimensional space of the Poncelet
polygons in the two-dimensional space: as is well-known, the
Poncelet polygons can be transformed by a projective transformation
into periodic trajectories of the billiards inside a quadric
in the two-dimensional space. Thus, the loop which started with Hitchin's observation in the elliptic case has been closed
for all genera.

The structure of the paper is as follows. In Section \ref{sect_Hurwitz}, we introduce the Hurwitz spaces of two-fold ramified covering over a Riemann sphere and basic tools for studying the variation of their branch points. In Section \ref{sect_Picard}, we discuss the Okamoto transformation taking the general solution of Painlev\'e's sixth equation with parameters $(0,0,0,\frac{1}{2})$ to the general solution of the equation with parameters \eqref{constants}.
In Section \ref{sect_EllipticOkamoto}, we introduce the main tool of our construction, the differential $\Omega$ on an elliptic curve, and show that it can be interpreted as the Okamoto transformation from Section \ref{sect_Picard}. We also show that the result of \cite{Hitchin} is a particular case of this transformation.
In Section \ref{sect_ESchlesinger}, we solve the rank two Schlesinger system associated with a family of elliptic curves in terms of the differential $\Omega$.
In Section \ref{sect_HSchlesinger}, we construct a differential $\Omega$ on a hyperelliptic curve as a generalization of the differential from Section \ref{sect_EllipticOkamoto}. Then we solve the rank two Schlesinger system associated with a family of hyperelliptic curves and discuss its relationship with the multidimensional Garnier system.
In Section \ref{sect_independence}, we show that the differential $\Omega$ is independent of the  choice of the canonical homology basis on the underlying curve and discuss the global behaviour of the constructed solutions.
In Section \ref{sect_tau}, we  compute the tau-function corresponding to our solution of the Schlesinger system from Section \ref{sect_ESchlesinger}, that is in the case of elliptic curves. We show that this tau-function coincides with the one computed in \cite{KiKo}.
In Section \ref{sect_backToPoncelet}, we discuss the relationship between our construction and the billiard ordered games.

\section{Hurwitz space}
\label{sect_Hurwitz}

Here we recall the definition of a Hurwitz space of hyperelliptic coverings of ${\mathbb  CP}^1$. A hyperelliptic covering is a pair $(\surf, f)$, where $\surf$ is a Riemann surface and $f:\surf \to {\mathbb  CP}^1$ is a function of degree two.
Zeros $P_1,\dots,P_n \in \surf$ of the differential $df$ are called ramification points and their images in ${\mathbb  CP}^1$ denoted here by $u_k=f(P_k)$ are called branch points of the covering.
 Two coverings $(\surf, f)$ and $(\tilde{\surf}, \tilde{f})$ are called equivalent if there exists a biholomorphic map $h: \surf \to \tilde{\surf}$ such that $f=\tilde{f}\circ h$. The Hurwitz space ${\mathcal H}_g^2$ is the space of equivalence classes of hyperelliptic ramified coverings. For degree two coverings, the equivalence amounts to the equality of the unordered sets of branch points, $\{u_1,\dots,u_n\}= \{\tilde{u}_1,\dots,\tilde{u}_n\} $. The set of branch points of a covering gives thus a set of coordinates on ${\mathcal H}_g^2$.
 The number of branch points of a genus $g$ hyperelliptic covering is given by the Riemann-Hurwitz formula: $n=2g+2$.

A Riemann surface for which a meromorphic function $f$ of degree two exists is called hyperelliptic. Such a surface $\mathcal L$ of genus $g$ can be seen as corresponding to a hyperelliptic algebraic curve (also denoted by $\mathcal L$) defined by the equation
\begin{equation*}
v^2=\prod_{i=1}^{2g+\epsilon}(u-u_i) \qquad \mbox{with}\;\;\;\; \epsilon\in\{1,2\}.
\end{equation*}
The corresponding hyperelliptic covering of $\mathbb{CP}^1$ is the pair $({\mathcal L}, f)$ where $f:{\mathcal L}\to \mathbb{CP}^1$ is given by $f((u,v))=u$. In what follows by abuse of notation we write $u( P)$ instead of $f(P )$ for a point $P=(u,v)$ on the curve.

Let $\{a_1,\dots,a_g; b_1,\dots,b_g\}$ be a canonical  homology basis on $\surf$, that is the intersection indices are as follows $a_j\circ a_k=0, \; b_j\circ b_k=0, \;a_k\circ b_j = \delta_{jk}$.
Denote by $\omega_1,\dots,\omega_g$ the holomorphic differentials on the surface $\surf$ normalized by the condition $\oint_{a_k}\omega_j=\delta_{jk}$. The $b$-periods of the differentials $\omega_k$ form the Riemann matrix ${\mathbb B}$ of the surface: ${\mathbb B}_{jk} = \oint_{b_k}\omega_j$.

Here we recall the Rauch formulas, see \cite{Fay92, KokoKoro}, describing the dependence of differentials on ${\mathcal L}$ on the variation of simple branch points of the covering. To this end, we introduce a bidifferential $W(P,Q)$ for $P,Q\in{\mathcal L}$, a differential with respect to $P$ and a differential with respect to $Q$, as the unique bidifferential on $\mathcal L$ having the following properties:
\begin{itemize}
\item Symmetry: $W(P,Q) = W(Q,P);$
\item No singularity except for a second order pole at the diagonal $P=Q$ with biresidue $1$: for $\xi$ being a local parameter near $P=Q$, the bidifferential has the following local expansion:
\begin{equation*}
W(P,Q) \underset{P\sim Q}{=} \left( \frac{1}{(\xi(P) - \xi(Q))^2}  + {\cal O}(1) \right)d\xi(P) d\xi(Q);
\end{equation*}
\item Vanishing of all $a$-periods: $\oint_{a_k} W(P,Q) = 0.$
\end{itemize}

As a consequence of this definition of $W$ we have: $\oint_{b_k}W(P,Q) = 2\pi{\rm i}\,\omega_k(P).$

Throughout the paper we use the following convention for the value of a differential at a given point. The evaluation is done with respect to the standard local parameters on a hyperelliptic curve: in a neighbourhood of a ramification point $P_j$, the standard local parameter is $\xi_j(P) = \sqrt{u(P) - u_j}$ if $u_j$ is finite and $\xi_\infty = \sqrt{u^{-1}}$ otherwise; in a neighbourhood of a regular point $Q$, the local parameter is $u$ if $u(Q)$ is finite and $u^{-1}$ otherwise.

Let us fix a point $R_0\in {\mathcal L}$ and let $\xi$ be the standard local parameter in a neighbourhood of $R_0$. An Abelian differential $\varphi(P)$ on ${\mathcal L}$ for $P\sim R_0$ has the form $\varphi(P) = g(\xi(P))d\xi(P)$ with some function $g(\xi)$. Then we put
\begin{equation}
\label{evaluation}
\varphi(R_0) := g(\xi(R_0)) = \frac{\varphi(P)}{d\xi(P)} \Big{|}_{P=R_0}
\end{equation}
and we say that {\it the differential $\varphi$ is equal to $\varphi(R_0)$} at the point $P=R_0$.

Any differential  $\varphi(P)$ on a hyperelliptic covering $({\mathcal L},u)$ of genus $g$ can be considered as an object defined in a neighbourhood of this covering in the Hurwitz space ${\mathcal H}^2_g$, and we write $\varphi = \varphi(P; u_1, \dots, u_n)$. When one of the coordinates $u_k$ varies, the complex structure on the surface ${\mathcal L}$ varies through the variation of the local parameter $\xi_k(P )=\sqrt{u(P )-u_k}$, which corresponds to moving locally in the Hurwitz space. The differential $\varphi$ varies accordingly. The Rauch variational formulas describe the dependence of the holomorphic normalized differentials $\omega_j(P)$, of $W(P,Q)$  and of the Riemann matrix ${\mathbb B}$ on such variations of the branch points $u_k$ provided the arguments $P,Q \in {\mathcal L}$ of the differential in question are fixed by the condition $u(P )=const, u(Q)=const$.

The Rauch variational formulas have the following form, see \cite{Fay92, KokoKoro}:
\begin{equation}
\label{RauchW}
\frac{\partial}{\partial u_k} W(P,Q)= \frac{1}{2} W(P,P_k) W(P_k,Q);
\end{equation}
\begin{equation}
\label{Rauchomega}
\frac{\partial \omega_j(P)}{\partial u_k} = \frac{1}{2} W(P,P_k) \omega_j(P_k);
\end{equation}
\begin{equation*}
\frac{\partial {\mathbb B}_{ij}}{\partial u_k} = \pi {\rm i}\, \omega_i(P_k) \omega_j(P_k),
\end{equation*}
where the evaluation of differentials at ramification points $P_k$ is performed  according to (\ref{evaluation}) with respect to the local parameter $\xi_k$.

\section{Picard solution and Okamoto transformation}
\label{sect_Picard}

Consider the transformed Weierstrass $\wp$-function with periods $2w_1$ and $2w_2$ satisfying the equation
\begin{equation}
\label{wp}
\left( {\wp^\prime}(z) \right)^2 = \wp(z)(\wp(z)-1)(\wp(z)-x).
\end{equation}
As is well known, see for example (3.6) in \cite{Okamoto} , the function
\begin{equation}
\label{Picard}
y_0 (x) = \wp(2c_1 w_1(x) + 2c_2 w_2(x))
\end{equation}
 is the Picard solution \cite{Picard} of Painlev\'e VI equation (\ref{Painleve}) with  constants $\alpha = \beta=\gamma=0,$ $\delta=1/2$.

The   transformation
\begin{equation}
\label{Okamoto}
y(x) = y_0 + \frac{y_0(y_0-1)(y_0-x)}{x(x-1)y_0^\prime - y_0(y_0-1)}
\end{equation}
which follows from \cite{Okamoto},  relates Picard's solution
$y_0(x)$ (\ref{Picard}) and the general solution $y(x)$ of the
Painlev\'e VI equation with constants (\ref{constants}).

Note that we obtained formula (\ref{Okamoto}) from Example 2.1 of \cite{Okamoto} after correcting the following misprint.  The expression for $\partial h/\partial t$ in Example 2.1 of \cite{Okamoto} should read:
\begin{equation}
\label{typo}
\frac{\partial h}{\partial t} = -q(q-1)p^2 + \{b_1(2q-1)-b_2 \}p-b_1^2.
\end{equation}
The very first $q$ in the right hand side of (\ref{typo}) is omitted in \cite{Okamoto}. In the notation of \cite{Okamoto} the Painlev\'e equations are written for a function $q(t)$, that is $q$ stands for what we denote by $y$ and the independent variable $t$ is the same as our $x$.

In the rest of this section we explain the situation in Example 2.1 of \cite{Okamoto} which is of interest to us and obtain the form of function (\ref{typo}) corresponding to this situation.

In \cite{Okamoto}, the set of parameters $(\alpha, \beta, \gamma, \delta)$ is transformed into two other sets of constants, $(\varkappa_0, \varkappa_1, \varkappa_\infty, \theta)$ and $(b_1, b_2, b_3, b_4)$, by the following rules:
\begin{equation*}
\alpha=\frac{1}{2}\varkappa_\infty^2, \qquad \beta=-\frac{1}{2}\varkappa_0^2, \qquad \gamma=\frac{1}{2}\varkappa_1^2, \qquad \delta = \frac{1}{2}(1-\theta^2);
\end{equation*}
and
\begin{equation*}
b_1=\frac{1}{2}(\varkappa_0+\varkappa_1), \qquad b_2=\frac{1}{2}(\varkappa_0-\varkappa_1), \qquad
b_3=\frac{1}{2}(\theta-1+\varkappa_\infty), \qquad b_4=\frac{1}{2}(\theta-1-\varkappa_\infty).
\end{equation*}
The transformation described in Example 2.1 of \cite{Okamoto}, let us denote it by ${\bf w}$, takes the general solution of the sixth Painlev\'e equation with parameters $(b_1, b_2, b_3, b_4)$ to the general solution of that equation with parameters $(b_3, b_2, b_1, b_4)$. The case interesting to us, the transformation of the Picard solution into the solution of the Painlev\'e VI with parameters (\ref{constants}), is a particular case of this transformation. Indeed, Picard's solution corresponds to  constants $\varkappa_0=\varkappa_1=\varkappa_\infty=\theta=0$ and $(b_1, b_2, b_3, b_4)=(0,0,-\frac{1}{2}, -\frac{1}{2})$. The transformation $\bf w$ maps this set of constants onto the set $(\tilde{b}_1, \tilde{b}_2,\tilde{b}_3,\tilde{b}_4)=(-\frac{1}{2}, 0, 0, -\frac{1}{2})$ and thus $\tilde{\varkappa}_0=-\frac{1}{2}, \tilde\varkappa_1=-\frac{1}{2}, \tilde\varkappa_\infty=\frac{1}{2}, \tilde\theta=\frac{1}{2}$, which gives the parameters from (\ref{constants}). Thus, in particular, ${\bf w}(y_0)=y$, where $y_0=y_0(x)$ is the Picard solution (\ref{Picard}) and $y=y(x)$ is the general solution of the Painlev\'e VI with parameters (\ref{constants}).

The Painlev\'e VI with constants $(0,0,0,\frac{1}{2})$ is considered in \cite{Okamoto} in its Hamiltonian formulation:
\begin{equation*}
\frac{dq}{dt} = \frac{\partial H}{\partial p}, \qquad \frac{dp}{dt} = -\frac{\partial H}{\partial q}
\end{equation*}
with the Hamiltonian
\begin{equation*}
H(t) = \frac{1}{t(t-1)} \left[  q(q-1)(q-t) p^2 +q(q-1)p +\frac{1}{4}(q-t) \right].
\end{equation*}
The transformation $\bf w$ is written with the help of an auxiliary function $h(t)$ defined by (1.6) of \cite{Okamoto} as follows:
\begin{equation*}
h(t) = t(t-1)H(t) +\sigma_2^\prime[b] t - \frac{1}{2} \sigma_2[b],
\end{equation*}
where $\sigma_2[b]$ is the elementary symmetric polynomial of degree $2$ in $b_1, b_2, b_3, b_4$, and $\sigma^\prime_2[b]$ is the elementary symmetric polynomial of degree $2$ in $b_1, b_3, b_4$. In our case, since $b_1=b_2=0,$ both polynomials coincide: $\sigma_2[b]=\sigma^\prime_2[b]= b_3b_4 = 1/4.$
Thus as is easy to see the $t$-derivative of the function $h(t)$ is given by (\ref{typo}) with $b_1=b_2=0$.

Now we can write the transformation $q_{\bf w}={\bf w}(q)$ given in Example 2.1 of \cite{Okamoto}, which in our notation corresponds to $y={\bf w}(y_0)$:
\begin{equation*}
q_{\bf w} = q+\frac{(b_3-b_1)q(q-1)p}{\partial h/\partial t + b_1^2}=q-\frac{1}{2p},
\end{equation*}
where the second equality is obtained by substituting the expression from (\ref{typo}) for the derivative of $h$ and the values of the parameters: $b_1=0$ and $b_3=-\frac{1}{2}$. Now we find $p$ in terms of the $t$-derivative of $q$ from the Hamiltonian system:
\begin{equation*}
2p=\left( \frac{t(t-1)\frac{dq}{dt}}{q(q-1)}-1\right)\frac{1}{q-t}.
\end{equation*}
Combining the last two equations and writing $y$ instead of $q_{\bf w}$ and $y_0$ instead of $q$ as well as $x$ instead of $t$, we get the Okamoto transformation in the form (\ref{Okamoto}).

\section{Elliptic curve and  Painlev\'e VI}
\label{sect_EllipticOkamoto}

The relationship between elliptic curves and the Painlev\'e VI
equations has a long history, from the classical work of
Picard \cite{Picard} and Fuchs \cite{Fuchs} to the contemporary papers of Hitchin \cite{Hitchin_Painleve, Hitchin} and Manin
\cite{Manin}. Some other aspects of this relationship have been treated in \cite
 {Guzzetti}.

We note that in order to use (\ref{Okamoto}) to write the solution $y(x)$ explicitly in terms of elliptic functions, we would need to express the derivative $dy_0/dx$ in elliptic functions. In this section, we find an expression for this derivative in terms of Abelian differentials on the associated elliptic curve, see (\ref{dy0_Omega}).

\subsection{Invariant form of the Okamoto transformation}

Let  $\surf=\surf(x)$ be the elliptic curve  associated with the Weierstrass $\wp$-function from Section \ref{sect_Picard}, that is the curve defined by equation
\begin{equation}
\label{ell_curve}
v^2=u(u-1)(u-x),
\end{equation}
where the coordinates are given by $(u,v)=(\wp(z), \wp^\prime(z))$.
We consider this curve as a two-fold ramified covering $({\mathcal L}, u)$ of the $u$-sphere ${\mathbb {CP}}^1$ and denote the ramification points of the covering as follows: $P_0:=(u=0, v=0)$, $P_1:=(u=1, v=0)$, $P_x:=(u=x, v=0)$ and $P_\infty= (\infty,\infty)$. The ramified covering can be presented in the form of a Hurwitz diagram where horizontal lines represent sheets of the covering and the vertical solid lines indicate ramification points.
Let $*$ stand, as before, for the elliptic involution interchanging the sheets of the covering so that for a point $P=(u,v)$ on the curve the result of the involution is $P^*=(u,-v)$.


 Let $a, b$ be a canonical homology basis on $\surf$. Denote by $\omega$ the holomorphic differential on $\surf$ normalized by the condition $\oint_a \omega=1$. Then its $b$-period is the period of the elliptic curve: $\oint_b\omega=\mu$ with  $\mu=2w_1/2w_2$. The elliptic curve $\surf$ is biholomorphic to its Jacobian $J(\surf):={\mathbb C}/\Lambda = {\mathbb C}/\{n+\mu m\mid n, m\in{\mathbb Z}\} $, the biholomorphism being given by the Abel map ${\mathcal A}:\surf\to J(\surf)$
\begin{equation*}
{\mathcal A}(P )=\int_{P_\infty}^P\omega,  \qquad P\in\mathcal L.
\end{equation*}
The transformed Weierstrass $\wp$-function (\ref{wp}) gives, after adjusting to the lattice generated by $2w_1$ and $2w_2$,  the inverse map: $\wp(2w_1{\mathcal A}(P )) = u(P )$ with $u(P )$ being the projection of the point $P\in\surf$ on the $u$-sphere.

Consider an arbitrary point
\begin{equation}
\label{z0}
z_0=c_1 + c_2\mu
\end{equation}
 in the Jacobian of the curve $\surf$ and let $Q_0$ be its preimage under the Abel map: ${\mathcal A}(Q_0)=z_0$, that is let the point $Q_0\in\surf$ be such that
\begin{equation}
\label{Q0}
\int_{P_\infty}^{Q_0}\omega = c_1 + c_2\mu.
\end{equation}
%

Note that the projection of $Q_0$ on the $u$-sphere is given by $y_0=\wp(2w_1z_0)$. Let us now consider the family of elliptic curves (\ref{ell_curve}) parametrized by the position of the branch point $x$. We define the point $z_0(x)$ in the Jacobian of the corresponding curve by equation (\ref{z0}) under the condition that the constants $c_1, c_2$ do not depend on $x$, that is $z_0(x)=c_1+c_2\mu(x).$ Recall that in this case the $u$-coordinate of $Q_0$, considered as a function of the moving branch point $x$, gives the Picard solution (\ref{Picard}) of a Painlev\'e VI equation.
The starting point of our construction is the point $Q_0$ and its elliptic involution, i.e., the two points $Q_0$ and $Q_0^*$  on $\surf$ whose $u$-coordinate coincides with the Picard solution: $u(Q_0)=u(Q_0^*) = y_0(x)$.

Note that we assume from now on that the constants $c_1, c_2$ are not equal to half-integers at the same time, that is $(c_1, c_2)\notin (\frac{1}{2}\mathbb Z)^2$, because otherwise our construction is impossible due to the fact the point $Q_0$ will coincide with one of the branch points and, therefore, with $Q_0^*$.

Let $\Omega_{Q_0,Q_0^*}$ be the normalized, $\oint_a \Omega_{Q_0,Q_0^*} =0$, differential of the third kind on $\surf$ with simple poles at $Q_0$ and $Q_0^*$ with the residues $1$ and $-1$, respectively. The following differential of the third kind on $\surf$ is the main tool in our construction:
\begin{equation}
\label{Omega}
\Omega = \Omega_{Q_0,Q_0^*} - 4 \pi {\rm i} c_2 \omega.
\end{equation}

The constants $c_1$ and $c_2$ define periods of $\Omega$ with respect to the canonical homology basis $a, b$ on the elliptic curve. Due to the normalization of $\Omega_{Q_0,Q_0^*}$ and $\omega$, the $a$-period of $\Omega$ is equal to $-4\pi{\rm i}c_2$ whereas the $b$-period is $4\pi{\rm i}c_1$. This can be seen as follows:
\begin{equation*}
\oint_{b}\Omega=2\pi{\rm i} \int_{Q_0^*}^{Q_0} \omega - 4 \pi {\rm i} c_2 \oint_b\omega= 4\pi{\rm i} \int_{P_\infty}^{Q_0} \omega - 4 \pi {\rm i} c_2 \mu,
\end{equation*}
where the first term is equal to $4\pi{\rm i} (c_1+c_2\mu)$ due to (\ref{Q0}). The first equality is a corollary of the Riemann bilinear relations and for the second equality we use $\omega(P^*)=-\omega(P )$ for any point $P$ on the elliptic curve.

In what follows, differentials are evaluated at specific points of the surface according to (\ref{evaluation}) with respect to the standard local parameters, namely: the parameter at $P_j$ is $\sqrt{u-j}$, with $j\in\{0,1,x\}$,  in a neighbourhood of $P_\infty$ the standard local parameter is $u^{-1/2}$ and in a neighbourhood of a regular point of the ramified covering the local parameter is $u$.

As an example of this evaluation, we compute  values of the holomorphic differential on $\mathcal L$ at various points of the curve. The quantities computed here will be used throughout the paper. The holomorphic differential on $\surf$ normalized by the condition $\oint_a \omega=1$ is expressed as follows in terms of the coordinates $z$ and $u$:
\begin{equation}
\label{holdiff}
\omega(P) = \frac{dz}{2w_1} =  \frac{du}{I_0\sqrt{u(u-1)(u-x)}}, \qquad P=(u,v) \in \surf, \;\; u = \wp(z), \;\;v=\wp'(z),
\end{equation}
with the normalization constant
\begin{equation}
\label{I0}
I_0:= \oint_a \frac{du}{\sqrt{u(u-1)(u-x)}} = \int_0^{2w_1} dz  = 2w_1 .
\end{equation}

The evaluation $\omega(P_j)$ of the holomorphic normalized differential (\ref{holdiff}) at ramification points $P_j$ with $j\in\{0,1,x\}$ with respect to the standard local parameters is as follows:
\begin{eqnarray}
&&\omega(P_0) := \frac{\omega(Q)}{d\sqrt{u(Q)}}{\Big |}_{Q=P_0, u=0} = \frac{2}{I_0\sqrt{x}}\,,
\label{omega_atP0}
\\
&&\omega(P_1) := \frac{\omega(Q)}{d\sqrt{u(Q)-1}}{\Big |}_{Q=P_1,u=1} = \frac{2}{I_0\sqrt{1-x}}\,,
\label{omega_atP1}
\\
&&\omega(P_x) := \frac{\omega(Q)}{d\sqrt{u(Q)-x}}{\Big |}_{Q=P_x, u=x} = \frac{2}{I_0\sqrt{x(x-1)}}.
\label{omega_atPx}
\end{eqnarray}
The evaluation of $\omega$ at the regular point $Q_0$ whose $u$-coordinate is $y_0$ is done with respect to the parameter $u$ and gives the following:
\begin{equation}
\label{omega_atQ0}
\omega(Q_0):= \frac{\omega(Q)}{du(Q)}{\Big |}_{Q=Q_0} =  \frac{1}{I_0\sqrt{y_0(y_0-1)(y_0-x)}}.
\end{equation}
Evaluation of other meromorphic differentials at given points is done similarly.

Now we are in a position to prove the main theorem of this Section which gives an invariant characterisation of the Okamoto transformation (\ref{Okamoto}) in terms of the associated elliptic curve.

\begin{theorem}
\label{thm_y}
Consider the family (\ref{ell_curve}) of elliptic curves $\surf(x)$ parametrized by $x$ and
let the constants $c_1, c_2\in \mathbb C$ be such that $(c_1, c_2)\notin (\frac{1}{2}\mathbb Z)^2$. Define a point $Q_0$ on a curve $\surf(x)$ from the family (\ref{ell_curve}) by condition \eqref{Q0} with $\mu$ being the period of $\surf(x)$ and let
 the differential $\Omega$ on $\surf(x)$ be defined by (\ref{Omega}). Then its zeros $P_y$ and $P_y^*$ project to the same point in the $u$-sphere and the function defined by $y(x)=u(P_y)=u(P_y^*)$ satisfies the Painlev\'e equation (\ref{Painleve}) with parameters (\ref{constants}).
\end{theorem}
{\it Proof.}
As is easy to see, the differential $\Omega$  in terms of the coordinate $u$ has the form:
%
\begin{equation}
\label{Omega_xy}
\Omega(P )= \frac{  \omega(P )}{\omega(Q_0)} \left[ \frac{1}{u(P )-y_0} - \frac{I}{I_0}  \right]- 4\pi{\rm i}c_2\omega(P ),
\end{equation}
where
\begin{equation}
I= \oint_a \frac{du}{(u-y_0) \sqrt{u(u-1)(u-x)}}.
\end{equation}
The quantity $I/I_0$ ensures the normalization of $\Omega_{Q_0Q_0^*}$ with respect to the $a$-cycle.

%
As is obvious from (\ref{holdiff}), the holomorphic differential $\omega$  has opposite signs on different sheets and thus from (\ref{Omega_xy}) we see that the same is true for our differential of the third kind: $\Omega(P^*) = - \Omega(P).$ Therefore its two zeros are mapped onto each other by the involution $*$ interchanging the sheets. We denote the zeros of $\Omega$ by $P_y$ and $P_y^*$, and their projection to the $u$-sphere, the base of the covering, by $y.$ The function $y(x)$ has the following form as can be found from (\ref{Omega_xy}):
\begin{equation}
\label{y}
\frac{1}{y-y_0} = \frac{I}{I_0} + 4\pi  {\rm i} c_2\omega(Q_0).
\end{equation}

We now prove that relations (\ref{Okamoto}) and (\ref{y}) between the functions $y$  and $y_0$ coincide. To this end, we find the derivative of $y_0$ with respect to the branch point $x$ as follows.

The Rauch variational formulas from Section \ref{sect_Hurwitz} take the following form for the covering $(\surf, u):$
\begin{equation}
\label{Rauchx}
\frac{d \omega(P)}{d x} = \frac{1}{2} W(P,P_x) \omega(P_x), \qquad \frac{d \mu}{d x} = \pi {\rm i}\, \omega^2(P_x).
\end{equation}
Using them to differentiate relation
(\ref{Q0}) with respect to the branch point $x$, we get:
\begin{equation*}
\omega(Q_0)y_0^\prime + \frac{1}{2} \int_{P_\infty}^{Q_0} W(P,P_x) \omega(P_x) = c_2 \pi {\rm i}\, \omega^2(P_x),
\end{equation*}
or equivalently
\begin{equation*}
\omega(Q_0)y_0^\prime + \frac{1}{4} \Omega_{Q_0Q_0^*}(P_x) \omega(P_x) = c_2 \pi {\rm i}\,  \omega^2(P_x),
\end{equation*}
where we used $2\int_{P_\infty}^{Q_0} W(P,P_x) = \int_{Q_0^*}^{Q_0} W(P,P_x)=\Omega_{Q_0Q_0^*}(P_x)$. This holds since $W(P,P_x) = - W(P^*,P_x)$ due to (\ref{WPPx}) and (\ref{holdiff}).
From this, using definition (\ref{Omega}) of the differential $\Omega$, we obtain for the derivative of $y_0(x)$:
\begin{equation}
\label{dy0_Omega}
\frac{dy_0}{dx} = -\frac{1}{4} \Omega(P_x) \frac{\omega(P_x)}{\omega(Q_0)}.
\end{equation}
Using (\ref{Omega_xy}) for $\Omega(P_x)$, we rewrite the derivative in the form:
\begin{equation}
\label{dy0}
\frac{d y_0}{d x} = \frac{1}{4} \frac{\omega^2(P_x)}{\omega^2(Q_0)} \left[4\pi {\rm i}\, c_2\,\omega(Q_0) - \frac{1}{x-y_0} + \frac{I}{I_0}\right].
\end{equation}
Using (\ref{dy0}) to express the right hand side of (\ref{y}) in terms of the $x$-derivative of $y_0$ and writing $\omega^2(P_x)$ and $\omega^2(Q_0)$ in terms of $x$ and $y_0$ as in (\ref{omega_atPx}), (\ref{omega_atQ0}), we arrive at (\ref{Okamoto}).
$\Box$

\subsection{Algebraic solutions of Painlev\'e VI}
\label{sect_Hitchin}

There was an intensive recent study of algebraic solutions to the
Painlev\'e VI equations. In \cite{DubrovinMazzocco} a
classification of  algebraic solutions to Painlev\'e VI with
parameters $(\alpha,0,0,\frac{1}{2})$
 for an arbitrary $\alpha$ was given. A more general classification  has been performed in \cite{LT}, see also \cite{Iwasaki} and references therein.
 The papers \cite{Manin}, \cite{Mazzocco}, \cite{DubrovinPainleve}, and \cite{DubrovinMazzocco1}  also discuss algebraic solutions of Painlev\'e VI.

 In this section we discuss the relationship between
Theorem \ref{thm_y} and Hitchin's work \cite{Hitchin} in which for
every closed Poncelet trajectory of length $k\in{\mathbb N}$,
$k\geq 3$ an algebraic solution to the Painlev\'e VI equation
(\ref{Painleve}), (\ref{constants}) is constructed in the
following way.

A closed Poncelet trajectory of length $k$ exists for two conics  defined by symmetric $3\times 3$ matrices $B$ and $C$ if and only if the point $(u,v)=(0,\sqrt{{\rm det}\,B})$ is of order $k$
on the elliptic curve  of the equation $v^2={\rm det}(B+uC)$, see \cite{GriffithsHarris}. That is if and only if $k\,{\mathcal A}((u,v))= \mu\, n+m$, with some integer constants $m$ and $n$, the Abel map $\mathcal A$ being based at the ramification point at infinity and $\mu$ being, as before, the period of the elliptic curve.

Given that one of the points, denote it by $\tilde{Q}_0$,  with $u$-coordinate equal to zero is of order $k$ on the elliptic curve, one can construct, as a corollary of the Abel theorem, a function $g(u,v)$ on the curve which has a zero of order $k$ at $\tilde{Q}_0$ and a pole of order $k$ at $u=\infty$ and no other zeros or poles. Then the function
\begin{equation*}
s(u,v)=\frac{g(u,v)}{g(u,-v)}
\end{equation*}
on the elliptic curve has a zero of order $k$ at $\tilde{Q}_0$ and a pole of order $k$ at $\tilde{Q}^*_0$ and no other zeros or poles.
The differential $ds$ has exactly two zeros away from the points $\tilde{Q}_0$ and $\tilde{Q}^*_0$. Since the involution $*$ takes
$s$ to $1/s$, these two zeros are paired by the involution, in other words, their $u$-coordinates coincide. Hitchin  proves \cite{Hitchin} that after a M\"obius transformation in the $u$-sphere taking the finite branch points of the curve to $0,1,x$ and preserving the point at infinity, the common $u$-coordinate of the two simple zeros of $ds$, as a function of $x$, defines an algebraic solution $y(x)$ to the Painlev\'e VI equation with constants (\ref{constants}).

In our terms, this situation after the M\"obius transformation is described as follows.

Let a point $Q_0$ on the curve (\ref{ell_curve}) be such that its image under the Abel map based at $P_\infty$ is a point of order $k\in {\mathbb N}$ in the Jacobian, that is
\begin{equation}
\label{kJacobian}
k\int_{P_\infty}^{Q_0}\omega= \mu n+m, \qquad \mbox{with some constant}\;\;\; m,n\in{\mathbb Z}.
\end{equation}
Then for $Q_0^*$ we have
\begin{equation*}
k\int_{P_\infty}^{Q_0^*}\omega= -\mu n-m.
\end{equation*}
Let $E(P,Q)$ be the prime form (see \cite{Fay92}) on the elliptic curve ${\mathcal L}$ (\ref{ell_curve}). Then the function $s$ with a zero of order $k$ at $Q_0$ and a pole of order $k$ at $Q_0^*$ and no other zeros or poles can be written as follows:
\begin{equation*}
s(P)  = \frac{E^k(P,Q_0)}{E^k(P,Q_0^*)} {\rm exp}\left\{-4n\pi {\rm i} \int_{Q_0^*}^P\omega\right\}   , \qquad P\in{\mathcal L.}
\end{equation*}
The differential $ds$ is therefore
\begin{equation*}
ds(P)  = ks(P )\left( d_P\log\frac{E(P,Q_0)}{E(P,Q_0^*)} - 4\pi{\rm i} \frac{n}{k}\omega(P ) \right) = ks(P ) \left( \Omega_{Q_0,Q_0^*}(P ) - 4\pi{\rm i} \frac{n}{k}\omega(P ) \right).
\end{equation*}
Thus according to \cite{Hitchin}, the projection $y(x)$ onto the $u$-sphere of two zeros of the differential
\begin{equation}
\label{differential}
\Omega_{Q_0,Q_0^*} - 4\pi{\rm i} \frac{n}{k}\omega
\end{equation}
on the elliptic curve ${\mathcal L}$ (\ref{ell_curve}) gives algebraic solutions to Painlev\'e equation (\ref{Painleve}), (\ref{constants}).

Note that if in our construction of Section \ref{sect_EllipticOkamoto} we put $c_1:=m/k$ and $c_2:=n/k$, then definition (\ref{Q0}) of the point $Q_0$ turns into (\ref{kJacobian}) and differential $\Omega$ (\ref{Omega}) becomes (\ref{differential}). In other words, Theorem \ref{thm_y} contains the result of \cite{Hitchin} outlined above as a special case of rational constants $c_1$ and $c_2$ from (\ref{z0}) provided the point $Q_0$ on the varying curve $\surf$ is defined by (\ref{kJacobian}) with constant $m, n\in \mathbb Z$.

\section{Elliptic curves and the Schlesinger system}

\label{sect_ESchlesinger}

In this Section we give a solution of the Schlesinger system (\ref{Schlesinger}) corresponding to the solution of the Painlev\'e equation (\ref{Painleve}), (\ref{constants}) from Theorem \ref{thm_y}. Note that the general solution to this Schlesinger system was already constructed in  \cite{Deift} and in \cite{KiKo} in terms of theta-functions. Here we give an alternative description of solutions to the four-point Schlesinger system in terms of differentials $\Omega$ (\ref{Omega}), $\omega$ (\ref{holdiff}) and the Picard solution $y_0$ (\ref{Picard}). We give an explicit and detailed proof using the Rauch variational formulas from Section \ref{sect_Hurwitz}.

Throughout this Section we assume that $y_0$ is the Picard solution (\ref{Picard}) with $(c_1, c_2)\notin (\frac{1}{2}\mathbb Z)^2$.

Note first that Theorem \ref{thm_y} shows that the only zero of the differential $\Omega$ (\ref{Omega}) in the $u$-sphere coincides with the only finite zero of the term $A_{12}(u)$ (\ref{A12}). Since the poles of $A_{12}(u)$ are known explicitly, we can construct this entry of the matrix $A(u)$ with the help of the differential $\Omega$ on the elliptic curve ${\mathcal L}$ (\ref{ell_curve}). Namely, the following proposition holds.

\begin{proposition}
\label{prop_A12}
Let $A^{(1)}, A^{(2)}, A^{(3)}\in sl(2,{\mathbb C})$ be matrices with eigenvalues $\pm 1/4$ satisfying the Schlesinger system (\ref{Schlesinger})  such that $A^{(1)}+A^{(2)}+A^{(3)}$ be a constant diagonal matrix. Then the $(12)$-term of the matrix $A$ (\ref{A}) is given by
\begin{equation}
\label{A12_omegas}
A_{12}(u) = \tilde{\kappa}(x) \frac{\Omega(P)}{\omega(P)} \;\frac{(u-y_0)}{u(u-1)(u-x)}, \qquad P\in {\mathcal L},
\end{equation}
where $\tilde{\kappa}(x)$ is a function of $x\in{\mathbb C}\setminus\{0,1,\infty\}$,  and $u=u(P)$ is the $u$-coordinate of the point $P\in{\mathcal L}$, with  ${\mathcal L}$ being the elliptic curve (\ref{ell_curve}).
\end{proposition}
\begin{remark}
{\rm We shall later see (Theorem \ref{thm_Schlesinger}) that  $\tilde{\kappa} = {1}/{I_0}$ in (\ref{A12_omegas}).}
\end{remark}
 {\it Proof.}  Note that, as can be deduced from  (\ref{Omega_xy}), the ratio
\begin{equation}
\label{Omega_omega}
\frac{\Omega(P)}{\omega(P)} =
\frac{1}{\omega(Q_0)} \left[ \frac{1}{u(P )-y_0} - \frac{I}{I_0}  \right] - { 4\pi{\rm i} c_2 }
\end{equation}
is a meromorphic function of $u\in{\mathbb {CP}^1}$, where $\omega(Q_0)$ is given by (\ref{omega_atQ0}), and $P$ is any finite point of the covering $({\mathcal L}, u)$. Therefore (\ref{A12_omegas}) as a function in the $u$-sphere has simple poles at $0,1$, and $x$, a simple zero at $u=y$ and a zero of order two at the point at infinity. We thus see that this function coincides with (\ref{A12}) up to a  factor independent of $u$.
$\Box$

\begin{corollary}
\label{cor_A12_residues}
The (12)-terms of residue-matrices $A^{(i)}$, $i=1,2,3$ are given by
\begin{equation}
\label{A12_residues_kappa}
A_{12}^{(1)} = -\frac{y_0}{x} \tilde{\kappa}\frac{\Omega(P_0)}{\omega(P_0)}, \qquad A_{12}^{(2)} = \frac{y_0-1}{x-1} \tilde{\kappa}\frac{\Omega(P_1)}{\omega(P_1)}, \qquad A_{12}^{(3)} = \frac{x-y_0}{x(x-1)}\tilde{\kappa} \frac{\Omega(P_x)}{\omega(P_x)}.
\end{equation}
They satisfy $A_{12}^{(1)}+A_{12}^{(2)}+A_{12}^{(3)} =0.$
\end{corollary}
{\it Proof.} These expressions are obtained by a straightforward calculation of residues of (\ref{A12_omegas}) at $u=0$, $u=1$ and $u=x$. Differential $A_{12}(u)du$ with $A_{12}$ given by (\ref{A12_omegas}) does not have a pole at the point at infinity, therefore (\ref{A12_residues_kappa}) gives all residues of this differential in the Riemann sphere and thus their sum is zero.
 $\Box$

In what follows, together with the normalized holomorphic differential $\omega$, we use the non-normalized holomorphic differential $\phi$ on the elliptic curve \eqref{ell_curve} defined by
\begin{equation}
\label{phi_ell}
\phi(P ) = I_0 \omega(P ) = \frac{du}{\sqrt{u \left( u - 1\right) \left( u - x\right)}} \qquad\mbox{with} \;\; P=(u, v).
\end{equation}
Note that the variation formula for $\phi$ can be easily obtained:
\begin{equation}
\label{phi_x}
\frac{d\phi(P )}{dx} = \frac{1}{2} \frac{\phi(P )}{u-x}
\end{equation}
where $u=u(P )$ is fixed under the variation of $x$ similarly to the Rauch formulas from Section~\ref{sect_Hurwitz}.

\begin{theorem}
\label{thm_Schlesinger}
 Consider the family  (\ref{ell_curve}) of elliptic curves with ramification points at $P_0,$ $P_1$, $P_x,$ and $P_\infty.$
Let, as before $y_0(x)$ be the Picard solution (\ref{Picard}) with
$(c_1, c_2)\notin (\frac{1}{2}\mathbb Z)^2$. For every value of $x$, let $Q_0$ and $Q_0^*$ be
the points on the corresponding elliptic curve having $y_0(x)$ as  $u$-coordinate:
$u(Q_0) = u(Q_0^*) = y_0(x)$, and let $\Omega(P)$ be the differential
(\ref{Omega}) of the third kind  with simple poles at $Q_0$ and
$Q_0^*.$ Let $\phi$ be the holomorphic non-normalized differential (\ref{phi_ell})
on the curve.

Introduce the following quantities:
\begin{equation}
\label{betas}
\beta_1:= -\frac{y_0}{x} \left( \frac{\Omega(P_0)}{\phi(P_0)} \right)^2, \qquad
\beta_2:= \frac{y_0-1}{x-1} \left( \frac{\Omega(P_1)}{\phi(P_1)} \right)^2, \qquad
\beta_3:= \frac{x-y_0}{x(x-1)} \left( \frac{\Omega(P_x)}{\phi(P_x)} \right)^2;
\end{equation}
and
\begin{equation}
\label{A12_residues}
A_{12}^{(1)} = -\frac{y_0}{x} \frac{\Omega(P_0)}{\phi(P_0)}, \qquad A_{12}^{(2)} = \frac{y_0-1}{x-1} \frac{\Omega(P_1)}{\phi(P_1)}, \qquad A_{12}^{(3)} = \frac{x-y_0}{x(x-1)} \frac{\Omega(P_x)}{\phi(P_x)}.
\end{equation}

The following matrices $A^{(i)}$ with $i=1,2,3$ give solution to the Schlesinger system (\ref{Schlesinger}):
\begin{equation}
\label{Ai}
A^{(i)}: = \left(
    \begin{array}{cc}
    -\frac{1}{4} -\frac{\beta_i}{2} & \qquad  A_{12}^{(i)}
    \\
    \\ - \frac{1}{4} \,\frac{\beta_i+\beta_i^2}{A_{12}^{(i)}} & \qquad \frac{1}{4} +\frac{\beta_i}{2}
    \end{array}\right).
\end{equation}
%
\end{theorem}

\begin{remark}
{\rm The eigenvalues of matrices $A^{(i)}$ (\ref{Ai}) are $\pm 1/4.$}
\end{remark}

\begin{remark}
\label{rmk_betas}
{\rm Note that using (\ref{omega_atP0})-(\ref{omega_atPx}) the coefficients $\beta_i$ (\ref{betas}) can also be rewritten in a simpler form:
\begin{equation*}
\beta_1:= -\frac{y_0}{4} \left( \Omega(P_0) \right)^2, \qquad
\beta_2:= \frac{1-y_0}{4} \left(\Omega(P_1) \right)^2, \qquad
\beta_3:= \frac{x-y_0}{4} \left( \Omega(P_x) \right)^2;
\end{equation*}
and the residues (\ref{A12_residues}) can be rewritten as follows:
\begin{equation}
\label{A12_residues_simpler}
A_{12}^{(1)} = -\frac{1}{4} y_0 \Omega(P_0)\phi(P_0), \qquad
A_{12}^{(2)} = \frac{1}{4} (1-y_0) \Omega(P_1)\phi(P_1), \qquad
A_{12}^{(3)} = \frac{1}{4} (x-y_0) \Omega(P_x)\phi(P_x).
\end{equation}
}
\end{remark}

In the proof of this theorem  we use a number of formulas grouped together in the next three lemmas.

\begin{lemma} Let ${\mathcal L}$ be the elliptic curve (\ref{ell_curve}), $P,Q$ denote regular points of the covering $({\mathcal L}, u)$, and $u(P)$ stand for the $u$-coordinate of the point $P$. Let the index $i$ stand for $0$, $1$ or $x$ and $u_i=i$ be the $u$-coordinate of ramification point $P_i$. Let $W$ be the bidifferential defined in Section \ref{sect_Hurwitz} and $\omega$ the holomorphic normalized differential (\ref{holdiff}).  Then the following formulas hold:
\label{lemma_formulas1}
\begin{equation}
\label{WPPx}
 W(P,P_i) = \left( \frac{1}{u(P )-u_i} - \frac{I^i}{I_0}\right) \frac{\omega(P)}{\omega(P_i)} ,
\end{equation}
where
\begin{equation}
\label{Ix}
I^i = \oint_a\frac{du}{(u-u_i) \sqrt{u(u-1)(u-x)}};
\end{equation}
and
\begin{equation}
\label{WPQPx}
W(P,P_i) \frac{\omega(Q)}{\omega(P)}= W(Q,P_i) + \left(\frac{1}{u(P) - u_i} - \frac{1}{u(Q) - u_i}\right) \frac{\omega(Q)}{\omega(P_i)}.
\end{equation}
%
\end{lemma}
{\it Proof.} Recall that $W(P,P_i)$ is a normalized ($\oint_aW(P,P_i) =0$) differential of the second kind with the only pole of second order at $P=P_i$ with biresidue equal to one. There is only one such differential, therefore it coincides with (\ref{WPPx}), (\ref{Ix}). Relation (\ref{WPQPx}) is a  corollary of (\ref{WPPx}), (\ref{Ix}).
$\Box$

\begin{lemma} Let ${\mathcal L}$ be the elliptic curve (\ref{ell_curve}), considered as a two-fold ramified covering of the $u$-sphere whose finite ramification points are $P_0$, $P_1$ and $P_x$.  Let $\Omega$ be differential (\ref{Omega}) of the third kind, and $\phi$ be the holomorphic non-normalized differential (\ref{phi_ell}). As before, the point $Q_0\in {\mathcal L}$ has coordinates $(y_0, \sqrt{y_0(y_0-1)(y_0-x)})$ with $y_0$ defined by (\ref{Picard}). Then the following formulas hold:
\label{lemma_formulas2}

\begin{eqnarray}
\label{Omega_0x}
&& \frac{\Omega(P_0)}{\phi(P_0)} - \frac{\Omega(P_x)}{\phi(P_x)} = \phi(Q_0) x(y_0 - 1);
\\ \nonumber
\\
%
\label{Omega_1x}
&& \frac{\Omega(P_1)}{\phi(P_1)} - \frac{\Omega(P_x)}{\phi(P_x)} = \phi(Q_0) y_0(x - 1);
\\ \nonumber
\\
\label{Omega_01}
&& \frac{\Omega(P_0)}{\phi(P_0)} - \frac{\Omega(P_1)}{\phi(P_1)} = \phi(Q_0) (y_0 - x).
\end{eqnarray}
\end{lemma}
{\it Proof.} Note that due to (\ref{Omega_omega}),
%
%
 the ratio
$\frac{\Omega(P)}{\phi(P)}$ is a function of $u(P )$ equal to $\frac{1}{\phi(Q_0)} \frac{1}{u(P)-y_0}$ plus the terms independent of $u(P)$.
Evaluating this ratio at ramification points and using (\ref{omega_atQ0}), we prove the lemma.
$\Box$

\begin{lemma} Let ${\mathcal L}$ be the elliptic curve (\ref{ell_curve}) and $P$  a regular point of the covering $({\mathcal L}, u)$. Let $\Omega$ be the differential (\ref{Omega}) of the third kind and $W$ the bidifferential defined in Section \ref{sect_Hurwitz}. The point $P_x$ is the ramification point with coordinates $(u=x,v=0)$ and $y_0$ given by (\ref{Picard}) is the $u$-coordinate of the points $Q_0$ and $Q_0^*$. Then the following variational formula holds:
\label{lemma_formulas3}
\begin{equation}
\label{dOmegax}
\frac{d}{dx} \Omega(P) = \frac{1}{2} \Omega(P_x) W(P,P_x) + \left(W(P,Q_0) - W(P,Q_0^*)\right) \frac{dy_0}{dx}.
\end{equation}
%
%
%
\end{lemma}
{\it Proof.}
Note that $\Omega$ can be written in the form:
\begin{equation*}
\Omega(P) = \int_{Q_0^*}^{Q_0} W(P,Q) - 4\pi {\rm i} \,c_2\omega(P).
\end{equation*}
We differentiate $W$ and $\omega$ according to the Rauch variational formulas (\ref{RauchW}) and (\ref{Rauchomega}). To differentiate the limits of integration, we note that the local parameter near $Q_0$ and $Q_0^*$ is $u$ and $u(Q_0)=u(Q_0^*)=y_0$. 

%
$\Box$

\begin{corollary} Let $u_i$ be the $u$-coordinate of the ramification point $P_i$ with $i\in\{ 0,1 \}$.
The following formulas hold for derivatives of $\Omega(P_i)$ :
\begin{equation}
\label{dOmegaPi}
\frac{d\Omega(P_i)}{dx} = \frac{1}{2} \frac{\Omega(P_x) \omega(P_x)}{\omega(P_i)} \frac{y_0-x}{(x-u_i)(y_0-u_i)}.
\end{equation}
\end{corollary}
{\it Proof.}
Note that since $P_i$ is a ramification point, the differential $W(P_i, P)$ has different  signs on different sheets and thus $W(P_i, Q_0)-W(P_i, Q_0^*)=2W(P_i, Q_0)$. This can be seen from expression  (\ref{WPPx}) for $W(P_i, P)$ in terms of the coordinate $u(P)$.  As is evident from (\ref{holdiff}), the holomorphic differential $\omega$ has different signs on different sheets, at points $P$ and $P^*$, therefore so does $W(P_i, P)$.

Now evaluating (\ref{dOmegax}) at $P=P_i$, we get
\begin{equation*}
\frac{d}{dx} \Omega(P_i) = \frac{1}{2} \Omega(P_x) W(P_i,P_x) + 2W(P_i,Q_0) \frac{dy_0}{dx}.
\end{equation*}
Using expression (\ref{dy0_Omega}) for the derivative of $y_0$, we rewrite this as follows
\begin{equation*}
\frac{d}{dx} \Omega(P_i) = \frac{1}{2} \Omega(P_x) \left( W(P_i,P_x) - W(P_i,Q_0) \frac{\omega(P_x)}{\omega(Q_0)}\right),
\end{equation*}
which coincides with (\ref{dOmegaPi}) due to (\ref{WPQPx}) of Lemma \ref{lemma_formulas1}.
$\Box$

{\it Proof of Theorem \ref{thm_Schlesinger}. }
We need to show that matrices (\ref{Ai}) satisfy differential equations (\ref{Schlesinger}). We do that by a straightforward differentiation of (\ref{Ai}) with respect to $x$.

 Let us start with the components $A_{12}^{(i)}$  of residue matrices $A^{(i)}$. The differential equations for them contained in (\ref{Schlesinger}) have the form:
 \begin{equation}
 \label{Schlesinger_12}
 \frac{d}{dx}A_{12}^{(1)} = \frac{2}{x}\left( A_{12}^{(1)} A_{11}^{(3)} - A_{12}^{(3)}A_{11}^{(1)}\right) \qquad \mbox{and} \qquad
  \frac{d}{dx}A_{12}^{(2)} = \frac{2}{x-1}\left( A_{12}^{(2)} A_{11}^{(3)} - A_{12}^{(3)}A_{11}^{(2)}\right),
 \end{equation}
where we took into account that the residue matrices are traceless. The differential equation for $A_{12}^{(3)}$ can be replaced by $A^{(1)}_{12}+A^{(2)}_{12}+A^{(3)}_{12} = 0$, which holds due to Corollary \ref{cor_A12_residues}.

We prove that $A_{12}^{(1)}$ from (\ref{A12_residues}) satisfies the first equation in (\ref{Schlesinger_12}). The proof for $A_{12}^{(2)}$ is entirely similar.


%

Let us differentiate expression (\ref{A12_residues}) for $A_{12}^{(1)}$ using the derivative (\ref{dOmegaPi}) of $\Omega(P_0)$, the variational formula (\ref{phi_x}) for $\phi(P )$, and the derivative (\ref{dy0_Omega}) of $y_0$:
\begin{multline*}
\frac{d}{dx}A_{12}^{(1)} =  -\frac{y_0}{x} \frac{1}{\phi(P_0)} \left(  \frac{1}{2}  \frac{\Omega(P_x)\omega(P_x)}{\omega(P_0)} \frac{(y_0-x)}{xy_0}   \right)
-\frac{ A_{12}^{(1)}}{x} +\frac{1}{2}\frac{ A_{12}^{(1)}}{x}
%
+\frac{1}{x^2(x-1)} \frac{\Omega(P_x)}{\phi{(P_x)}}\frac{\Omega(P_0)}{\phi{(P_0)}} \frac{1}{\phi(Q_0)},
\end{multline*}
where (\ref{omega_atPx}) was also used.
Then using again (\ref{omega_atP0}), (\ref{omega_atPx}), and the form (\ref{A12_residues}) for  $A_{12}^{(3)}$, we arrive at:
\begin{equation*}
\frac{d}{dx}A_{12}^{(1)} =  \frac{ A_{12}^{(3)}}{2x} - \frac{ A_{12}^{(1)}}{2x}   +\frac{1}{x^2(x-1)} \frac{\Omega(P_x)}{\phi{(P_x)}}\frac{\Omega(P_0)}{\phi{(P_0)}} \frac{1}{\phi(Q_0)}.
\end{equation*}
The last term in the right hand side can be rewritten using first (\ref{omega_atQ0}) and then (\ref{Omega_0x}) as follows:
\begin{equation*}
\frac{\phi(Q_0)y_0(y_0-1)(y_0-x)}{x^2(x-1)} \frac{\Omega(P_x)}{\phi{(P_x)}}\frac{\Omega(P_0)}{\phi{(P_0)}}  = \frac{y_0(y_0-x)}{x^3(x-1)} \frac{\Omega(P_x)}{\phi{(P_x)}}\frac{\Omega(P_0)}{\phi{(P_0)}}  \left(   \frac{\Omega(P_0)}{\phi{(P_0)}} -  \frac{\Omega(P_x)}{\phi{(P_x)}} \right),
\end{equation*}
which is equal to $\frac{A_{12}^{(3)}}{x} \beta_1 - \frac{A_{12}^{(1)}}{x}\beta_3.$ This gives
\begin{equation}
\label{dA_12}
\frac{d}{dx}A_{12}^{(1)} =  \frac{ A_{12}^{(3)}}{x} \left( \frac{1}{2} + \beta_1\right) - \frac{ A_{12}^{(1)}}{x} \left( \frac{1}{2} + \beta_3\right),
\end{equation}
which coincides with (\ref{Schlesinger_12}).

\vskip 1cm

Let us now verify the differential equations for the $(11)$- and $(22)$-entries of the matrices $A^{(i)}$. Recall that $A_{11}^{(i)} = - A_{22}^{(i)}$. These equations are of the form:
 \begin{equation}
 \label{Schlesinger_11}
 \frac{d}{dx}A_{11}^{(1)} = \frac{1}{x}\left( A_{12}^{(3)} A_{21}^{(1)} - A_{12}^{(1)}A_{21}^{(3)}\right) \qquad \mbox{and} \qquad
 \frac{d}{dx}A_{11}^{(2)} = \frac{1}{x}\left( A_{12}^{(3)} A_{21}^{(2)} - A_{12}^{(2)}A_{21}^{(3)}\right).
 \end{equation}
The  equation for $A_{11}^{(3)}$ can be replaced by $A^{(1)}_{11}+A^{(2)}_{11}+A^{(3)}_{11} = -\frac{1}{4}$ which can be proven by a straightforward calculation using formulas (\ref{Omega_0x}) - (\ref{Omega_01}).

As before, we prove the first equation in (\ref{Schlesinger_11}), the second equation is verified analogously. Rewriting (\ref{Schlesinger_11}) in terms of $\beta_i$, we get
 \begin{equation}
 \label{dbeta1_A}
 \frac{d\beta_1}{dx}= \frac{1}{2x} \left(  \frac{A_{12}^{(3)}}{A_{12}^{(1)}} (\beta_1+\beta_1^2) - \frac{A_{12}^{(1)}}{A_{12}^{(3)}} (\beta_3+\beta_3^2)  \right);
 \end{equation}
the same equation in terms of differentials $\Omega$ and $\phi$ becomes:
 \begin{equation}
  \label{dbeta1}
 \frac{d\beta_1}{dx}= \frac{1}{2x} \frac{\Omega(P_0)}{\phi(P_0)} \frac{\Omega(P_x)}{\phi(P_x)}   \left(  \frac{x-y_0}{x(x-1)} \left\{1-\frac{y_0}{x} \left( \frac{\Omega(P_0)}{\phi(P_0)} \right)^2\right\} + \frac{y_0}{x} \left\{1+ \frac{x-y_0}{x(x-1)} \left( \frac{\Omega(P_x)}{\phi(P_x)} \right)^2 \right\}  \right).
 \end{equation}
Now we need to obtain (\ref{dbeta1}) by differentiating $\beta_1$ defined by (\ref{betas}). Note that $\beta_1$ can be written as follows:
 \begin{equation}
 \label{beta1_A12}
\beta_1 = -\frac{x}{y_0} \left( A_{12}^{(1)} \right)^2 ,
 \end{equation}
therefore we can use (\ref{dA_12}) for the derivative of $A_{12}^{(1)}$ when differentiating $\beta_1$. Thus we get
 \begin{equation}
 \label{dbeta1_Rauch}
 \frac{d\beta_1}{dx}= -\frac{1}{y_0} \left( A_{12}^{(1)} \right)^2 -\frac{x}{4y_0^2} \left( A_{12}^{(1)} \right)^2 \Omega(P_x)\frac{\omega(P_x)}{\omega(Q_0)} - 2A_{12}^{(1)} \frac{1}{y_0} \left\{       { A_{12}^{(3)}} \left( \frac{1}{2} + \beta_1\right) - { A_{12}^{(1)}} \left( \frac{1}{2} + \beta_3\right)         \right \},
 \end{equation}
where (\ref{dy0_Omega}) was used for $dy_0/dx,$ which can be rewritten by using (\ref{omega_atPx}) and (\ref{omega_atQ0}) as follows:
 \begin{equation*}
 \frac{dy_0}{dx} = -\frac{1}{4} \Omega(P_x) \frac{\omega(P_x)}{\omega(Q_0)} = -\frac{1}{x(x-1)} \frac{\Omega(P_x)}{\phi(P_x)} \frac{1}{\phi(Q_0)}.
 \end{equation*}

Proving that (\ref{dbeta1_Rauch}) coincides with (\ref{dbeta1}) is a lengthy but straightforward computation. One can, for example, rewrite both expressions only in terms of the quantity $\frac{\Omega(P_0)}{\phi(P_0)}$  using formulas from Lemma \ref{lemma_formulas2} and see, with the help of (\ref{omega_atQ0}),  that terms with different powers of this quantity in both expressions coincide.

\vskip 1cm

It remains to show that Schlesinger system (\ref{Schlesinger}) holds also for the $(21)$-components of matrices $A^{(i)}$ (\ref{Ai}), that is
\begin{equation}
\label{Schlesinger_21}
\frac{dA_{21}^{(1)}}{dx} = \frac{2}{x} \left(  A_{21}^{(3)} A_{11}^{(1)} - A_{21}^{(1)} A_{11}^{(3)} \right) \qquad \mbox{and} \qquad
\frac{dA_{21}^{(2)}}{dx} = \frac{2}{x-1} \left(  A_{21}^{(3)} A_{11}^{(2)} - A_{21}^{(2)} A_{11}^{(3)} \right),
\end{equation}
 the  differential equation for $A_{12}^{(3)}$ being replaced by $A_{21}^{(1)}+A_{21}^{(2)}+A_{21}^{(3)}=0.$

Equations (\ref{Schlesinger_21}) can be derived from the form of matrices (\ref{Ai}) without even using the Rauch variational formulas. Namely, we find the derivative of $\beta_1$ from ${d\beta_1}/{dx} = -2 {dA_{11}^{(1)}}/{dx}$ and (\ref{Schlesinger_11}) for the derivative of $A_{11}^{(1)}$. In this way we have
\begin{equation*}
\frac{dA_{21}^{(1)}}{dx} = \frac{1}{xA_{12}^{(1)}} \left(     A_{12}^{(3)} A_{21}^{(1)} - A_{12}^{(1)}A_{21}^{(3)}  \right) \left(\frac{1}{2} + \beta_1\right) - \frac{A_{21}^{(1)}}{A_{12}^{(1)}} \frac{dA_{12}^{(1)}}{dx}.
\end{equation*}
Now replacing the derivative of  $A_{12}^{(1)}$ from (\ref{dA_12}) and simplifying, we arrive at (\ref{Schlesinger_21}) for $A_{21}^{(1)}$. The differential equation for $A_{21}^{(2)}$ is proven analogously.

To prove that the sum of $(21)$-terms vanishes, we rewrite it as follows:
\begin{equation*}
-4(  A_{21}^{(1)} +A_{21}^{(2)} + A_{21}^{(3)}  ) = \frac{\Omega(P_0)}{\phi(P_0)} (1+\beta_1) + \frac{\Omega(P_1)}{\phi(P_1)} (1+\beta_2) + \frac{\Omega(P_x)}{\phi(P_x)} (1+\beta_3)
\end{equation*}
and show that it is equal to zero by using formulas from Lemma \ref{lemma_formulas2} and (\ref{omega_atQ0}).

This completes the proof of Theorem \ref{thm_Schlesinger}.
$\Box$

\begin{corollary}
For solution (\ref{Ai}) of the Schlesinger system, the residue of the matrix $A\,du$ at the point $u=\infty$ is
\begin{equation*}
A^{(\infty)} = - A^{(1)} - A^{(2)} - A^{(3)} = \left(\begin{array}{cc}\frac{1}{4} & 0 \\ \\0 & -\frac{1}{4}\end{array}\right).
\end{equation*}
\end{corollary}
{\it Proof.} This follows from the proof of Theorem \ref{thm_Schlesinger}, where we proved $A_{12}^{(1)}+A_{12}^{(2)}+A_{12}^{(3)}=0,$ $A^{(1)}_{11}+A^{(2)}_{11}+A^{(3)}_{11} = -\frac{1}{4}$, and $A_{21}^{(1)}+A_{21}^{(2)}+A_{21}^{(3)}=0.$   Recall that ${\rm tr} A^{(i)}=0$.
$\Box$

\section{Hyperelliptic curves and the Schlesinger system}
\label{sect_HSchlesinger}

\subsection{Schlesinger system}
Consider now the linear matrix system (\ref{linsys_intro})
for a $2\times 2$ matrix function $\Phi(u)$ defined in the Riemann sphere, with the matrix $A\in sl(2,{\mathbb C})$ given by (\ref{Ah}) having simple poles at $2g+1$ points $u_1, \dots, u_{2g+1}$ in the complex plane with residues  $A^{(i)}\in sl(2,{\mathbb C})$ and a simple pole at the point at infinity.
%
Isomonodromic deformations of this system are described by Schlesinger system (\ref{hSchlesinger_intro}). 
%
As in the four-point Schlesinger system from Section \ref{sect_ESchlesinger}, the eigenvalues $\frac{t_i}{2}$ and $-\frac{t_i}{2}$ of the residue-matrices $A^{(i)}$ are integrals of motion of the Schlesinger system. Again,  one can consider solutions to system (\ref{hSchlesinger_intro}) up to a conjugation by a constant invertible matrix. This freedom is used to assume that the matrix $A^{(\infty)}:=-(A^{(1)} + \dots + A^{(2g+1)})$ is diagonal, $A^{(\infty)}={\rm diag}(\lambda, -\lambda)$.

In what follows, we verify explicitly that the functions we construct solve Schlesinger differential equations (\ref{hSchlesinger_intro}). For our proof it is convenient to see these equations written in the components as follows:
\begin{equation}
\label{Schlesinger_matrix}
\frac{\partial A^{(j)}}{\partial u_k} = \frac{1}{u_k-u_j}\left(\begin{array}{cc}A_{12}^{(k)} A_{21}^{(j)} - A_{12}^{(j)} A_{21}^{(k)} & 2( A_{12} ^{(j)} A_{11}^{(k)} - A_{12}^{(k)} A_{11}^{(j)} ) \\
\\
2( A_{21} ^{(k)} A_{11}^{(j)} - A_{21}^{(j)} A_{11}^{(k)} ) & A_{12}^{(j)} A_{21}^{(k)} - A_{12}^{(k)} A_{21}^{(j)}\end{array}\right) .
\end{equation}

\subsection{Jacobi inversion on a hyperelliptic curve}
\label{subsection_inversion}
Schlesinger system (\ref{hSchlesinger_intro})  is naturally associated with  the family of  hyperelliptic curves
\begin{equation}
\label{hyperelliptic}
v^2= \prod_{j=1}^{2g+1} (u-u_j).
\end{equation}
Our main goal in this Section is to  find a solution to the Schlesinger system in terms of quantities defined on these curves by generalizing the approach developed in Section \ref{sect_ESchlesinger} for elliptic curves.

For a point $P$ on a hyperelliptic curve $\surf$ from the family (\ref{hyperelliptic}), we write $P=(u,v)$. As before, we represent the curve as a two-fold ramified covering $(\surf, u)$ of the $u$-sphere and denote by $P_k$ the ramification points of the curve: $P_k=(u_k,0).$

Suppose a canonical homology basis $\{a_k, b_k\}_{k=1}^g$ is chosen on a hyperelliptic curve from family (\ref{hyperelliptic}) and  consider an arbitrary point $z_0$ in the Jacobian of the curve:
\begin{equation}
\label{preinversion}
z_0: = c_1 + {\mathbb B}c_2,
\end{equation}
where $\mathbb B$ is the corresponding Riemann matrix. The Jacobi
inversion of this point is given by a positive divisor $D =
Q_1+\dots + Q_g$ of degree $g$ on the curve:
\begin{equation}
\label{inversion}
\int_{P_\infty}^{Q_1} {\bf \omega} + \dots + \int_{P_\infty}^{Q_g} {\bf \omega} = c_1+{\mathbb B}c_2,
\end{equation}
where ${\bf \omega}$ is a vector of normalized holomorphic differentials: ${\bf \omega} = (\omega_1,\dots,\omega_g)^t$ with $\oint_{a_k}\omega_j=\delta_{kj}.$

Let us suppose that all points $Q_j$ are different from ramification points.
We introduce the following analogue of  the differential of the third kind $\Omega$ (\ref{Omega}) for the hyperelliptic curve:
\begin{equation}
\label{hOmega}
\Omega(P ) = \sum_{j=1}^g \Omega_{Q_j Q_j^*}(P ) - 4\pi {\rm i}\, c_2^t{\bf \omega}(P ).
\end{equation}
Here $\Omega_{Q_j Q_j^*}$ is the normalized differential of the third kind with zero $a$-periods and with poles at $Q_j$ and $Q_j^*$ of residues $1$ and $-1$, respectively. This differential can be written in terms of the bidifferential $W$ (see Section \ref{sect_Hurwitz}):
\begin{equation}
\label{OmegaW}
\Omega_{Q_j Q_j^*}( P) = \int_{Q^*_j}^{Q_j} W(P,R).
\end{equation}
Note that analogously to the genus one case, the $a$- and $b$-periods of $\Omega$ are given by the vectors of constants $c_1$ and $c_2$. Clearly, the $a_k$-period is equal to $-4\pi {\rm i}c_{2k}$ where $c_{2k}$ is the $k$th component of $c_2$. For the $b$-periods we have
\begin{equation*}
(\oint_{b_1}\Omega, \dots, \oint_{b_g}\Omega)^t = 2\pi{\rm i}  \sum_{j=1}^g \int_{Q_j^*}^{Q_j} \omega - 4\pi {\rm i}\, {\mathbb B} c_2 = 4\pi{\rm i}  \sum_{j=1}^g \int_{P_\infty}^{Q_j} \omega - 4\pi {\rm i}\, {\mathbb B} c_2=4\pi{\rm i} c_1.
\end{equation*}
Here the first equality is a corollary of the Riemann bilinear relations, in the second equality we used the fact that the holomorphic differentials change sign under the hyperelliptic involution, $\omega(P^*)=-\omega(P )$, where the involution is defined for a point $P=(u,v)$ of the curve by $(u,v)^*=(u,-v).$ The last equality is due to (\ref{inversion}).



Denote by $q_j$ the $u$-coordinate of the point $Q_j$, $j=1, \dots, g$. We now consider a family of the hyperelliptic curves (\ref{hyperelliptic}), that is we allow the branch points vary. As before, the Jacobian of the curve varies accordingly and we define the point $z_0$ in the varying Jacobian by (\ref{preinversion}) with the assumption that the vectors $c_1$ and $c_2$ are kept fixed, that is independent of the branch points $\{u_k\}_{k=1}^{2g+1}$. In other words, the point $z_0$ will move in the complex plane in such a way that its coordinates with respect to the lattice $\Lambda$ of the Jacobian stay constant under the variation of the curve.

In this paper we suppose  that the divisor $D$ is such that
\begin{itemize}
\item[(i)] all points $Q_j$ are distinct,
\item[(ii)] none of them coincides with a ramification point of the surface, and
\item[(iii)] no two points $Q_i$ and $Q_k$ are paired by the hyperelliptic involution, that is $Q_i\neq Q_k^*$ for any $i,k=1,\dots,g$.
\end{itemize}
The family of hyperelliptic curves (\ref{hyperelliptic}) forms a fiber bundle over the base space $\Pi$ of unordered sets of $2g+1$ distinct points $\{u_1, \dots, u_{2g+1}\}$. Suppose our curve $\mathcal L$ from the family (\ref{hyperelliptic}) and the point $z_0$ (\ref{preinversion}) in its Jacobian are such that the corresponding divisor $D$ satisfies requirements (i)-(iii). Denote by ${\mathcal U}=\{u_1, \dots, u_{2g+1}\}$ the projection of the fiber $\surf$ on the base of the bundle. Then, for the fixed coordinates $c_1,$ $c_2$ of the point $z_0$,  there exists an open neighbourhood of  $\mathcal U$ in $\Pi$ such that for any point $\tilde{\mathcal U}$ in this neighbourhood, the corresponding fiber $\tilde{\mathcal L}$ is a surface for which the Jacobi inversion of the point $z_0: = c_1 + \tilde{\mathbb B}c_2$ also satisfies conditions (i)-(iii). In this paper, we consider small deformations of $\mathcal U$ which keep it in the described neighbourhood.

We now want to obtain an expression for derivatives of $q_j$ with respect to branch points $u_k$  generalizing (\ref{dy0_Omega}) to the case of hyperelliptic curves.  Differentiating (\ref{inversion}) with respect to a branch point $u_k$, using the Rauch variational formulas (\ref{Rauchomega}), we get:
\begin{equation*}
\frac{1}{2} {\bf \omega}(P_k) \sum_{j=1}^g \int_{P_\infty}^{Q_j} W(P,P_k)  + \sum_{j=1}^g{\bf \omega}(Q_j)\frac{\partial q_j}{\partial u_k} = \pi {\rm i} \, {\bf \omega} (P_k) {\bf \omega}^t (P_k) c_2.
\end{equation*}
As before, the differential $W(P,P_k)$ changes sign under the hyperelliptic involution: $W(P^*,P_k) = -W(P,P_k)$, and therefore $\int_{Q^*_j}^{Q_j} W(P,P_k)=2\int_{P_\infty}^{Q_j} W(P,P_k)$.
Thus, due to (\ref{OmegaW}), we can rewrite the last equation as follows:
\begin{equation*}
 \sum_{j=1}^g{\bf \omega}(Q_j)\frac{\partial q_j}{\partial u_k} = - \frac{1}{4} {\bf \omega} (P_k) \left( \sum_{j=1}^g \Omega_{Q_jQ_j^*}(P_k) - 4 \pi {\rm i} \, {\bf \omega}^t (P_k) c_2\right) = -\frac{1}{4} {\bf \omega}(P_k) \Omega(P_k),
\end{equation*}
with $\Omega$ given by (\ref{hOmega}). This implies
\begin{equation}
\label{der}
 \sum_{j=1}^g\tilde{\omega}(Q_j)\frac{\partial q_j}{\partial u_k} = -\frac{1}{4} \tilde{\omega} (P_k) \Omega(P_k),
\end{equation}
for   any holomorphic differential $\tilde{\omega}$ on the hyperelliptic curve.

Denote by $\phi$ the following non-normalized holomorphic differential on the hyperelliptic curve:
\begin{equation}
\label{phi}
\phi = \frac{du}{ \sqrt{\prod_{j=1}^{2g+1}(u-u_j)}}.
\end{equation}

\begin{remark}
{\rm Note that  the evaluation of $\phi$ at a ramification point $P_n$ with respect to the standard local parameter $\sqrt{u-u_n}$ gives
%
$\phi(P_n) = 2 \left(\sqrt{\prod_{j=1, j\neq n}^{2g+1}(u_n-u_j)}\right)^{-1}$
%
and its derivatives with respect to branch points are
\begin{equation}
\label{phi_derivatives}
\frac{\partial \phi(P_n)}{\partial u_k} = \frac{1}{2} \frac{\phi(P_n)}{u_n-u_k}, \qquad k\neq n.
\end{equation}
}
\end{remark}

The differentials $\phi_0=\phi$ and  $\phi_l = u^l\phi$ with $1\leq l\leq g-1$ are holomorphic on the surface.
%
%
Since relation (\ref{der}) is valid for $\tilde{\omega} = \phi_l$, by linearity we get
\begin{equation}
\label{der_poly}
 \sum_{j=1}^g{\mathbb P}(q_j) \phi(Q_j)\frac{\partial q_j}{\partial u_k} = -\frac{{\mathbb P}(u_k)}{4} \phi (P_k) \Omega(P_k),
\end{equation}
where ${\mathbb P}$ is a polynomial of degree at most $g-1$.

In the space of holomorphic differentials on the hyperelliptic surface $\mathcal L$,  in addition to the basis of normalized differentials $\omega_1, \dots, \omega_g$ let us consider the following basis:
\begin{equation}
\label{vj}
v_j( P) = \frac{\phi( P) \prod^g_{\alpha=1, \alpha\neq j}(u-q_\alpha)}{\phi(Q_j)\prod^g_{\alpha=1, \alpha\neq j}(q_j-q_\alpha)}, \qquad j=1,\dots, g
\end{equation}
where $P \in \surf$ is a point on the surface and $u$ its projection on the base of the covering, the $u$-sphere. Differential $v_j$ vanishes at $2g-2$ points $Q_i, Q^*_i$ with $i=1,\dots,g$, $i\neq j$ and its evaluation with respect to the standard local parameter $u$ at the point $Q_j$ according to (\ref{evaluation}) gives one: $v_j(Q_j)=1$. Clearly, the differentials $v_j$ change sign under the hyperelliptic involution: $v_j(P^*) = - v_j(P ). $ Note that any holomorphic differential $\tilde{\omega}$ can be decomposed with respect to this basis as follows:
\begin{equation}
\label{v-decomposition}
\tilde{\omega}(P )= \sum_{j=1}^g \tilde{\omega}(Q_j) v_j(P ).
\end{equation}
To see this, it suffices to consider the difference of the left and right hand sides of (\ref{v-decomposition}). This  difference is a holomorphic differential vanishing at $2g$ points $Q_i, Q^*_i$ with $i=1,\dots,g$ and therefore is zero.

Note that for $g=1$ this new basis reduces to the differential $v_1 = \frac{\omega(P )}{\omega(Q_0)}$ with $\omega$ being the normalized holomorphic differential (\ref{holdiff}). Thus we get the following hyperelliptic analogue of (\ref{dy0_Omega}).

\begin{lemma}
\label{lemma_derivatives}
Derivatives of the $u$-coordinates of the points of the divisor $D = Q_1+\dots + Q_g$  (\ref{inversion}) with respect to the branch points are given by
\begin{equation}
\label{derivatives}
\frac{\partial q_j}{\partial u_k} = -\frac{1}{4}  \Omega(P_k) v_j(P_k).
\end{equation}
\end{lemma}
{\it Proof.}
This follows immediately from (\ref{der}) with $\tilde{\omega}=v_j$ since $v_j$ vanishes at all the points of the divisor $D$ except $Q_j$, where it is equal  to one.
$\Box$


\subsection{Solution of the Schlesinger system}

The following function is an analogue  in the hyperelliptic case of $A_{12}$ given by (\ref{A12}) with $\kappa=1$, or of $A_{12}$ given by (\ref{A12_omegas}) with $\tilde{\kappa} = 1/I_0$:
\begin{equation}
\label{A12_hyperelliptic}
A_{12}(u) =  \frac{\Omega(P )}{\phi(P )} \frac{\prod_{\alpha=1}^g(u-q_\alpha)}{\prod_{j=1}^{2g+1}(u-u_j)}.
\end{equation}
The differential $A_{12}(u)du$ on the $u$-sphere has zeros at the zeros of $\Omega$ and simple poles at the branch points of the curve as its only singularities. Its residues at the simple poles should be the $12$-entries of the residue matrices $A^{(j)}$:
\begin{equation}
\label{def_residues}
A_{12}^{(n)} = \frac{1}{4} \Omega(P_n) \phi(P_n) \prod_{\alpha=1}^g (u_n-q_\alpha).
\end{equation}
Note that these residues provide a straightforward generalization of the elliptic ones (\ref{A12_residues_simpler}).

\begin{remark}
\label{rmk_sum12}{\rm
Note that the sum or matrices (\ref{def_residues}) vanishes: $\sum_{n=1}^{2g+1}A^{(n)}_{12}=0$ as a sum of residues of the differential $A_{12}(u)du$.}
\end{remark}

Let us now introduce the following quantities for $n=1, \dots, 2g+1$:
\begin{equation}
\label{betas_hyperelliptic}
\beta_n:=
\frac{1}{4} \, \Omega(P_n) \sum_{j=1}^g v_j(P_n) - \frac{1}{2}\Omega(P_\infty)A_{12}^{(n)},
\end{equation}
where $\Omega(P_\infty)$ is the evaluation of the differential $\Omega$ at the ramification point at infinity with respect to the local parameter $u^{-1/2}$, and $A_{12}^{(n)}$ is given by (\ref{def_residues}).

\begin{theorem}
\label{thm_main}
Consider the family (\ref{hyperelliptic}) of hyperelliptic curves with variable branch points $\{u_k\}_{k=1}^{2g+1}$.
Let $\beta_n$ and $A_{12}^{(n)}$ be defined respectively by (\ref{betas_hyperelliptic}) and  (\ref{def_residues}) for this family of curves.
 Then the following matrices $A^{(n)}$ as functions of the branch points  $\{u_k\}_{k=1}^{2g+1}$ solve the Schlesinger system (\ref{hSchlesinger_intro}):
\begin{equation}
\label{Ai_hyper}
A^{(n)}: = \left(
    \begin{array}{cc}
    -\frac{1}{4} -\frac{\beta_n}{2} & \qquad  A_{12}^{(n)}
    \\
    \\ - \frac{1}{4} \,\frac{\beta_n+\beta_n^2}{A_{12}^{(n)}} & \qquad \frac{1}{4} +\frac{\beta_n}{2}
    \end{array}\right).
\end{equation}
\end{theorem}

\begin{remark}\rm{ As in the elliptic case, the eigenvalues of the residue matrices are $\pm1/4$, or in other words,
${\rm tr} \left(A^{(n)}\right)^2=1/8$.
 This can be seen directly from the form of matrices $A^{(n)}$.}
\end{remark}

\begin{remark}
\rm{
Here we would like to discuss definition (\ref{betas_hyperelliptic}) for coefficients $\beta_n$. We show that $\beta_n$ from (\ref{betas_hyperelliptic}) for a surface of genus one turn into (\ref{betas}). And thus, as is easy to see, (\ref{Ai_hyper})  in genus one becomes the solution  (\ref{Ai}) of the Schlesinger system in the elliptic case.

First, note that the differential $\Omega$ (\ref{hOmega}) can be written using the basis (\ref{vj}) of holomorphic differentials on the surface as follows:
\begin{equation}
\label{hOmega_vj}
\Omega(P ) = \sum_{j=1}^gv_j(P )\left(\frac{1}{u-q_j}   + \alpha_j\right)
\end{equation}
with some normalizing coefficients $\alpha_j \in \mathbb C$.

Now, multiplying $j$th term in the sum by $(u-q_j)$ and evaluating at $P=P_n, u=u_n$ we get the quantity which in genus one turns into $\Omega(P_n)(u_n-y_0)$:
\begin{equation}
\label{temp_betas}
 \sum_{j=1}^gv_j(P_n )\left(1   + \alpha_j(u_n-q_j)\right).
\end{equation}
The coefficients $\beta_n$ (\ref{betas_hyperelliptic}) are obtained by multiplying (\ref{temp_betas}) by $\Omega(P_n)/4$, which in genus one gives the multiplication of $\Omega(P_n)(u_n-y_0)$ by $\Omega(P_n)/4$:
\begin{equation}
\label{temp_betas_def}
\beta_n=\frac{1}{4}\Omega(P_n) \sum_{j=1}^gv_j(P_n )\left(1   + \alpha_j(u_n-q_j)\right).
\end{equation}
To see that this coincides with the definition of $\beta_n$, we note that \eqref{temp_betas_def} can  be rewritten as
\begin{equation*}
\beta_n=\frac{1}{4}\Omega(P_n) \sum_{j=1}^gv_j(P_n )  + A_{12}^{(n)} \sum_{j=1}^g
 \frac{\alpha_j}{\phi(Q_j)\prod^g_{\alpha=1, \alpha\neq j}(q_j-q_\alpha)}.
\end{equation*}
Evaluating (\ref{hOmega_vj}) at the ramification point at infinity, which we denote by $P_\infty$, we get
\begin{equation}
\label{Omega_infty}
\Omega(P_\infty) = \sum_{j=1}^g\alpha_jv_j(P_\infty ) = -2\sum_{j=1}^g\frac{ \alpha_j}{\phi(Q_j)\prod^g_{\alpha=1, \alpha\neq j}(q_j-q_\alpha)}
\end{equation}
and thus prove that (\ref{temp_betas_def}) coincides with (\ref{betas_hyperelliptic}). From this discussion we get immediately that in genus one $\beta_n$ (\ref{temp_betas_def}) turns into $\beta_n=\frac{1}{4}\Omega^2(P_n)(u_n-y_0)$, which, according to Remark \ref{rmk_betas}, coincides with (\ref{betas}).
}
\end{remark}

Here is another interesting property of coefficients $\beta_n$.

\begin{lemma}
\label{lemma_betasum}
 For a hyperelliptic surface of genus $g$
\begin{equation*}
\sum_{n=1}^{2g+1} \beta_n = -g.
\end{equation*}
\end{lemma}
{\it Proof.} Given definition (\ref{betas_hyperelliptic}) for $\beta_n$, we have
\begin{equation*}
\sum_{n=1}^{2g+1} \beta_n = \sum_{j=1}^g  \sum_{n=1}^{2g+1} \frac{ A_{12}^{(n)} }{ \phi(Q_j)(u_n-q_j)\prod^g_{\alpha=1, \alpha\neq j}(q_j-q_\alpha) } - \frac{1}{2}\Omega(P_\infty)\sum_{n=1}^{2g+1} A_{12}^{(n)}.
\end{equation*}
Note that $\sum_{n=1}^{2g+1} A_{12}^{(n)}=0$ as sum of residues of the differential $A_{12}(u)du$ where $A_{12}(u)$ is given by (\ref{A12_hyperelliptic}).

Let us consider a slightly modified differential on the base of the covering: $\frac{A_{12}(u)}{u-q_j}du$. This differential vanishes at the point at infinity and has simple poles at $u=q_j$ and at all finite branch points $u_k$, $k=1,\dots,2g+1$ of the hyperelliptic surface with the following residues:
\begin{equation*}
\underset{u=u_n}{\rm res}\frac{A_{12}(u)}{u-q_j}du = \frac{A_{12}^{(n)}}{u_n-q_j}; \qquad
\underset{u=q_j}{\rm res}\frac{A_{12}(u)}{u-q_j}du = \frac{\prod_{\alpha=1, \alpha\neq j}^g(q_j-q_\alpha)}{\phi(Q_j ) \prod_{m=1}^{2g+1}(q_j-u_m)} = \phi(Q_j) \prod_{\alpha=1, \alpha\neq j}^g(q_j-q_\alpha),
\end{equation*}
where we used that $\phi^{2}(Q_j) = \prod_{m=1}^{2g+1}(q_j-u_m)^{-1}$.

Since the sum of residues vanishes, we get
\begin{equation*}
\sum_{n=1}^{2g+1} \frac{A_{12}^{(n)}}{u_n-q_j} =  - \phi(Q_j) \prod_{\alpha=1, \alpha\neq j}^g(q_j-q_\alpha),
\end{equation*}
and therefore
\begin{equation*}
\sum_{n=1}^{2g+1} \frac{A_{12}^{(n)}}{(u_n-q_j)\phi(Q_j) \prod_{\alpha=1, \alpha\neq j}^g(q_j-q_\alpha)} =   -1,
\end{equation*}
which completes the proof.
$\Box$

\begin{corollary}
\label{cor_sum11}
The following relation holds
\begin{equation*}
\sum_{i=1}^{2g+1} A_{11}^{(i)} = -\sum_{i=1}^{2g+1} A_{22}^{(i)} = -\frac{1}{4}.
\end{equation*}
\end{corollary}
{\it Proof.} This follows immediately from Lemma \ref{lemma_betasum} and definition (\ref{Ai_hyper}) of matrix coefficients $A_{11}^{(i)}$, $A_{22}^{(i)}$.
$\Box$

\subsection{Proof of Theorem \ref{thm_main}}

The next two lemmas, Lemmas \ref{lemma_prod} and
\ref{lemma_Omega_derivative}, contain technical results which will
be used in the proof of Theorem \ref{thm_main}. Lemmas
\ref{lemma_diffeq_12} - \ref{lemma_sum21}  and Corollary
\ref{cor_Ainfty} contain our proof of Theorem \ref{thm_main}.

\begin{lemma}
\label{lemma_prod}
For every $k\in\{1, \dots, 2g+1\}$  the following equality holds
\begin{equation*}
\frac{1}{\prod_{\alpha=1}^g(u_k-q_\alpha)} = \sum_{j=1}^g \frac{1}{(u_k-q_j)\prod_{\alpha=1, \alpha\neq j}^g(q_j-q_\alpha)}.
\end{equation*}
\end{lemma}
{\it Proof.}  Decomposing the holomorphic differential $\phi$ with respect to the basis $\{v_j\}$ as in (\ref{v-decomposition}), we get
\begin{equation*}
\phi(P_k) = \sum_{j=1}^g \phi(Q_j) v_j(P_k).
\end{equation*}
Substituting now the expression for $v_j$ from (\ref{vj})  we obtain:
\begin{equation*}
\phi(P_k) = \sum_{j=1}^g \frac{\phi(P_k)\;\prod_{\alpha=1}^g(u_k-q_\alpha)}{(u_k-q_j) \prod_{\alpha=1, \alpha\neq j}^g(q_j-q_\alpha)}.
\end{equation*}
Dividing both sides by $\phi(P_k)\;\prod_{\alpha=1}^g(u_k-q_\alpha)$ we prove the lemma.
$\Box$

\begin{corollary}
\label{corrolary_prod}
The following equality holds for any distinct $n, k \in \{1, \dots, 2g+1\}$
\begin{equation}
\label{prod}
\sum_{j=1}^g \phi(Q_j) \frac{\partial q_j}{ \partial u_k} \prod_{\alpha=1, \alpha\neq j}^g (u_n-q_\alpha) = - \frac{1}{4} \frac{\Omega(P_k) \phi(P_k)}{u_n-u_k} \left(  \prod_{\alpha=1}^g (u_n-q_\alpha) - \prod_{\alpha=1}^g (u_k-q_\alpha)  \right).
\end{equation}
\end{corollary}
{\it Proof.}  Let us slightly rewrite the left hand side of (\ref{prod}) and then use (\ref{derivatives}) and (\ref{vj}) to substitute an algebraic expression for the derivatives ${\partial q_j}/{ \partial u_k}$:
\begin{equation*}
\prod_{\alpha=1}^g(u_n-q_\alpha) \sum_{j=1}^g \phi(Q_j) \frac{\partial q_j}{ \partial u_k} \frac{1}{u_n-q_j}
=
- \frac{1}{4} \prod_{\alpha=1}^g(u_n-q_\alpha)  \sum_{j=1}^g \frac{\phi(P_k)\Omega(P_k)\prod_{\alpha=1}^g(u_k-q_\alpha)}{(u_k-q_j) (u_n-q_j) \prod_{\alpha=1, \alpha\neq j}^g(q_j-q_\alpha)}.
 \end{equation*}

Now using Lemma \ref{lemma_prod}, we note that the difference of products in the right hand side of (\ref{prod}) can be rewritten as
\begin{equation*}
  \prod_{\alpha=1}^g (u_n-q_\alpha) - \prod_{\alpha=1}^g (u_k-q_\alpha)
=
\prod_{\alpha=1}^g (u_n-q_\alpha)\prod_{\alpha=1}^g (u_k-q_\alpha)
 \sum_{j=1}^g \frac{(u_n-u_k)}{(u_k-q_j) (u_n-q_j) \prod_{\alpha=1, \alpha\neq j}^g(q_j-q_\alpha)}.
\end{equation*}
Comparing the right hand sides of the last two equalities we prove the statement.
$\Box$

\begin{lemma}
\label{lemma_Omega_derivative}
The following formulas hold for derivatives of the differential $\Omega$ (\ref{hOmega}) with respect to the branch points:
\begin{equation}
\label{Omega_derivative}
\frac{\partial \Omega(P_n)}{\partial u_k} = \frac{1}{2}\frac{\Omega(P_k)}{u_k-u_n} \frac{\phi(P_k)}{\phi(P_n)}\frac{\prod_{\alpha=1}^g (u_k-q_\alpha)}{\prod_{\alpha=1}^g (u_n-q_\alpha)};
\end{equation}
\begin{equation}
\label{Omega_infty_derivative}
\frac{\partial \Omega(P_\infty)}{\partial u_k} = -\frac{1}{4}{\Omega(P_k)}   \phi(P_k)    \prod_{\alpha=1}^g (u_k-q_\alpha).
\end{equation}
\end{lemma}
{\it Proof.}
By definition (\ref{hOmega}) of the differential $\Omega$ and because of (\ref{OmegaW}), using the Rauch variational formulas (\ref{RauchW}), we get
\begin{equation}
\label{OmegaR_derivative}
\frac{\partial \Omega( R)}{\partial u_k} = \frac{1}{2} W(R,P_k) \Omega(P_k) + \sum_{\alpha=1}^g 2W(R,Q_\alpha) \frac{\partial q_\alpha}{\partial u_k}.
\end{equation}
We are going to use this equality with $R=P_n$ or $R=P_\infty.$
Thus we need the following representation of the differentials $W(P,P_n)$ and $W(P,P_\infty)$ in terms of the coordinate $u$ on the base of the covering:
\begin{equation}
\label{hW}
W(P, P_n) = \frac{1}{u(P ) -u_n } \frac{\phi(P )}{\phi(P_n)} - \sum_{j = 1}^g \frac{I_{u_n}^{a_j}}{\phi(P_n)} \omega_j(P )
\end{equation}
where $I_{u_n}^{a_j}$ is the following integral over the cycle $a_j$:
\begin{equation*}
I_{u_n}^{a_j} = \oint_{a_j} \frac{du}{(u-u_n)  \sqrt{\prod_{\alpha=1}^{2g+1}(u-u_\alpha)}}
\end{equation*}
and
\begin{equation}
\label{hW_infty}
W(P, P_\infty) = -\frac{u^g\phi(P )}{2} + \sum_{j = 1}^g \left(\oint_{aj}\frac{u^g\phi(P )}{2}\right) \omega_j(P ).
\end{equation}
Evaluating (\ref{hW}) at $P=P_k$ and $P=Q_\alpha$, we rewrite the derivative (\ref{OmegaR_derivative}) of $\Omega(P_n)$ in the form:
\begin{multline*}
\frac{\partial \Omega(P_n)}{\partial u_k}= \frac{1}{2}\Omega(P_k) \left(  \frac{1}{u_k-u_n} \frac{\phi(P_k)}{\phi(P_n)} - \sum_{j=1}^g  \frac{I_{u_n}^{a_j}}{\phi(P_n)} \omega_j(P_k)\right)\\
+2\sum_{\alpha=1}^g \left(  \frac{1}{q_\alpha-u_n} \frac{\phi(Q_\alpha)}{\phi(P_n)} - \sum_{j=1}^g \frac{I^{a_j}_{u_n}}{\phi(P_n)} \omega_j(Q_\alpha)  \right) \frac{\partial q_\alpha}{\partial u_k}.
\end{multline*}
Changing the order of summation in the last term we single out the sum $\sum_{\alpha=1}^g \omega_j(Q_\alpha) \frac{\partial q_\alpha}{\partial u_k}$ which we replace due to (\ref{der}) with $-\omega_j(P_k) \Omega(P_k)/4$ and obtain after some simplifications:
\begin{equation*}
\frac{\partial \Omega(P_n)}{\partial u_k}= \frac{1}{2} \frac{\Omega(P_k)}{u_k-u_n} \frac{\phi(P_k)}{\phi(P_n)}
+\frac{2}{\phi(P_n)}\sum_{\alpha=1}^g  \frac{\phi(Q_\alpha)}{q_\alpha-u_n}  \frac{\partial q_\alpha}{\partial u_k}.
\end{equation*}
Now multiplying and dividing the last term in the right hand side by $\prod_{j=1}^g (u_n-q_j)$ and using (\ref{prod}), we prove (\ref{Omega_derivative}).

Doing a similar calculation for derivative (\ref{OmegaR_derivative}) with $R=P_\infty$, we arrive at the following:
\begin{equation*}
\frac{\partial \Omega( P_\infty)}{\partial u_k} = -\frac{1}{4}u_k^g\Omega(P_k)\phi(P_k) - \sum_{j=1}^g q_j^g\phi(Q_j)\frac{\partial q_j}{\partial u_k}.
\end{equation*}
Due to (\ref{der_poly}) the first term in the right hand side can be replaced by $u_k\sum_{j=1}^g q_j^{g-1}\phi(Q_j) \frac{\partial q_j}{\partial u_k}$, and thus combining terms and using (\ref{derivatives}) and (\ref{vj}) for the derivatives of $q_j$, we get
\begin{equation*}
\frac{\partial \Omega( P_\infty)}{\partial u_k} = -\frac{1}{4}{\Omega(P_k)}   \phi(P_k)    \prod_{\alpha=1}^g (u_k-q_\alpha) \sum_{j=1}^g \frac{q_j^{g-1}}{\prod_{\alpha=1, \alpha\neq j}^g(q_j-q_\alpha)}.
\end{equation*}
The sum in the right hand side is equal to one as can be seen from equating to zero the sum of residues of the form
\begin{equation*}
\frac{u^{g-1}du}{\prod_{\alpha=1}^g (u-q_\alpha)}
\end{equation*}
on the Riemann sphere.
$\Box$

Now we begin to prove Theorem \ref{thm_main}.
\begin{lemma}
\label{lemma_diffeq_12}
The residues $A_{12}^{(i)}$ satisfy the following differential equation:
\begin{equation}
\label{diffeq_12_hyperelliptic}
\frac{\partial A_{12}^{(n)}}{\partial u_k}  =   \frac{A_{12}^{(k)}}{u_k-u_n} \left(\beta_n + \frac{1}{2}\right)
-\frac{A_{12}^{(n)}}{u_k-u_n} \left(\beta_k + \frac{1}{2}\right) .
\end{equation}
\end{lemma}

{\it Proof.}
This differential equation is just the $12$-component of equation (\ref{Schlesinger_matrix}) for the matrices $A^{(i)}$ defined by (\ref{Ai_hyper}).
Note that the derivative (\ref{Omega_derivative}) of $\Omega(P_n)$  can be rewritten in the form
\begin{equation*}
\frac{\partial \Omega(P_n)}{\partial u_k} = \frac{2A_{12}^{(k)}}{u_k-u_n} \frac{1}{\phi(P_n)\prod_{\alpha=1}^g (u_n-q_\alpha)}.
\end{equation*}
Using this and the derivatives (\ref{phi_derivatives}) of $\phi(P_n)$, a straightforward differentiation of $A_{12}^{(n)}$ (\ref{def_residues}) with respect to $u_k$ gives:
\begin{equation}
\label{temp}
\frac{\partial A_{12}^{(n)}}{\partial u_k}
= \frac{1}{2} \frac{A_{12}^{(k)}}{u_k-u_n} - \frac{1}{2}\frac{A_{12}^{(n)}}{u_k-u_n}
- \frac{1}{4} \Omega(P_n) \phi(P_n) \prod_{\alpha=1}^g(u_n-q_\alpha) \sum_{\alpha=1}^g \frac{1}{u_n-q_\alpha} \frac{\partial q_\alpha}{\partial u_k}.
\end{equation}
Note also that using (\ref{derivatives}) and (\ref{vj}) for  $\frac{\partial q_\alpha}{\partial u_k}$ the last term in (\ref{temp}) can be written as
\begin{multline}
\label{temp_A12}
-\frac{1}{4} \Omega(P_n) \phi(P_n) \prod_{\alpha=1}^g(u_n-q_\alpha) \sum_{\alpha=1}^g \frac{1}{u_n-q_\alpha} \frac{\partial q_\alpha}{\partial u_k}
= \sum_{\alpha=1}^g \frac{A_{12}^{(n)}A_{12}^{(k)}}{\phi(Q_\alpha)(u_n-q_\alpha)(u_k-q_\alpha)\prod_{j=1,j\neq \alpha}^g(q_\alpha-q_j)}.
\end{multline}

For $\beta_n$ defined by (\ref{betas_hyperelliptic}), we have
\begin{equation*}
\frac{\beta_n}{A_{12}^{(n)}} =   \sum_{j=1}^g\frac{1}{\phi(Q_j)(u_n-q_j)\prod^g_{\alpha=1, \alpha\neq j}(q_j-q_\alpha)}  - \frac{1}{2}\Omega(P_\infty),
\end{equation*}
and therefore
\begin{equation}
\label{boston}
\frac{\beta_n}{A_{12}^{(n)}}-\frac{\beta_k}{A_{12}^{(k)}}= \sum_{\alpha=1}^g \frac{u_k-u_n}{\phi(Q_\alpha)(u_n-q_\alpha)(u_k-q_\alpha)\prod_{j=1,j\neq \alpha}^g(q_\alpha-q_j)}.
\end{equation}
This shows that the last term in (\ref{temp}) is equal to
\begin{equation*}
\frac{A_{12}^{(n)}A_{12}^{(k)}}{u_k-u_n}\left( \frac{\beta_n}{A_{12}^{(n)}}-\frac{\beta_k}{A_{12}^{(k)}} \right),
\end{equation*}
and thus (\ref{temp}) and (\ref{diffeq_12_hyperelliptic}) coincide.
$\Box$

\begin{remark}
{\rm Note that equality (\ref{boston}) is a hyperelliptic analogue of relations from Lemma \ref{lemma_formulas2}. }
\end{remark}

\begin{lemma}
\label{lemma_diffeq_11}
The diagonal entries $A_{11}^{(i)}$ of the residue matrices $A^{(i)}$ (\ref{Ai_hyper}) satisfy the following differential equation:
\begin{equation*}
\frac{\partial A_{11}^{(n)}}{\partial u_k} = \frac{1}{u_k-u_n} \left(  A_{12}^{(k)} A_{21}^{(n)} -  A_{12}^{(n)} A_{21}^{(k)}  \right)
\end{equation*}
where the $(21)$-entries $A_{21}^{(i)}$  of the matrices $A^{(i)}$ are also defined by (\ref{Ai_hyper}).
\end{lemma}

{\it Proof.} This differential equation is equivalent the following equation for $\beta_n:$
\begin{equation}
\label{beta_diffeq}
\frac{\partial \beta_n}{\partial u_k} =
\frac{1}{2}\frac{1}{u_k-u_n} \left(  (\beta_n+\beta_n^2) \frac{A_{12}^{(k)}}{A_{12}^{(n)}} - (\beta_k+\beta_k^2) \frac{A_{12}^{(n)}}{A_{12}^{(k)}}  \right).
\end{equation}

Now we are going to differentiate expressions (\ref{betas_hyperelliptic}) for $\beta_n$ with respect to branch points $u_k$ of the hyperelliptic curve in a straightforward way using (\ref{diffeq_12_hyperelliptic}) for derivatives of $A_{12}^{(n)}$, (\ref{phi_derivatives}) for derivatives of $\phi(P_n)$, Lemma \ref{lemma_Omega_derivative} for derivatives of $\Omega(P_n)$ and $\Omega(P_\infty)$ and (\ref{derivatives}) for derivatives of $q_j$. We have to show that the result of this differentiation coincides with the right hand side of (\ref{beta_diffeq}). The proof is  technically involved and thus we include it with a lot of detail.

By  writing $\beta_n=A_{12}^{(n)} \frac{\beta_n}{A_{12}^{(n)}}$ and using (\ref{diffeq_12_hyperelliptic}) to differentiate $A_{12}^{(n)}$ we have
\begin{multline}
\label{beta_n_temp}
\frac{\partial \beta_n}{\partial u_k} =
\left(
   \frac{A_{12}^{(k)}}{u_k-u_n} \left(\beta_n + \frac{1}{2}\right)
-\frac{A_{12}^{(n)}}{u_k-u_n} \left(\beta_k + \frac{1}{2}\right)
\right) \frac{\beta_n}{A_{12}^{(n)}}
\\
+ A_{12}^{(n)} \frac{\partial}{\partial u_k} \left(
\sum_{j=1}^g \frac{1}{\phi(Q_j)(u_n-q_j)\prod_{\alpha=1,\alpha\neq j}^g(q_j-q_\alpha)} - \frac{\Omega(P_\infty)}{2}
\right).
\end{multline}
Note that (\ref{Omega_infty_derivative}) implies that $\partial_{u_k} \Omega(P_\infty) = -A_{12}^{(k)}$. Using definition (\ref{phi}) of $\phi$ and (\ref{derivatives})  for  $\partial_{u_k}q_\alpha$ we get
\begin{equation*}
\frac{\partial \phi^{-1}(Q_j)}{\partial u_k}  = \frac{1}{2\phi(Q_j)} \left(  \frac{1}{u_k-q_j} - \frac{1}{4} \sum_{l=1}^{2g+1} \frac{\Omega(P_k)}{q_j-u_l}v_j(P_k)  \right).
\end{equation*}
Similarly we differentiate the product of $(q_j-q_\alpha)$:
\begin{equation*}
\frac{\partial}{\partial u_k}  \prod_{\alpha=1,\alpha\neq j}^g(q_j-q_\alpha) = -\frac{1}{4} \prod_{\alpha=1,\alpha\neq j}^g(q_j-q_\alpha) \Omega(P_k) \sum_{l=1, l\neq j}^g \frac{v_j(P_k) - v_l(P_k)}{q_j-q_l}.
\end{equation*}
Thus, using (\ref{boston}) and definition (\ref{vj}) of differentials $v_j$, the derivative of $\beta_n$ becomes
\begin{multline*}
\frac{\partial \beta_n}{\partial u_k} =
\frac{1}{u_k-u_n}\left(
   \frac{A_{12}^{(k)}}{A_{12}^{(n)}} \left(\beta^2_n + \frac{\beta_n}{2}\right)
-\beta_n\beta_k - \frac{\beta_n}{2}\right)
+ \frac{1}{2}\frac{A_{12}^{(n)}}{u_k-u_n} \left( \frac{\beta_n}{A_{12}^{(n)}} - \frac{\beta_k}{A_{12}^{(k)}}  \right)
\\
- \frac{\Omega(P_k)\Omega(P_n)}{16} \sum_{j=1}^g  v_j(P_n) v_j(P_k)
\left(  \frac{1}{2}\sum_{l=1}^{2g+1} \frac{1}{q_j-u_l}
 +   \frac{1}{u_n-q_j} - \sum_{m=1, m\neq j}^g\frac{1}{q_j-q_m}
 \right)
 \\
 -\frac{\Omega(P_k)\Omega(P_n)}{16} \sum_{j=1}^g  v_j(P_n)\sum_{l=1, l\neq j} \frac{v_l(P_k)}{q_j-q_l} + \frac{1}{2}A_{12}^{(k)}A_{12}^{(n)}.
\end{multline*}

Now in the term $\beta_n\beta_k$  we replace $\beta_k$ due to (\ref{boston}) by
\begin{equation*}
\beta_k = \frac{A_{12}^{(k)}}{A_{12}^{(n)}}\beta_n    -    \frac{u_k-u_n}{4}\Omega(P_n)\frac{A_{12}^{(k)}}{A_{12}^{(n)}}
\sum_{j=1}^g \frac{v_j(P_n)}{u_k-q_j},
\end{equation*}
which leads to the following simplification in the first line for $\partial_{u_k}\beta_n$
\begin{multline*}
\frac{\partial \beta_n}{\partial u_k} =
\frac{1}{2}\frac{1}{u_k-u_n}\left(
   \frac{A_{12}^{(k)}}{A_{12}^{(n)}} \beta_n - \frac{A_{12}^{(n)}}{A_{12}^{(k)}} \beta_k   \right)
+ \frac{\Omega(P_n)}{4}A_{12}^{(k)}\frac{\beta_n}{A_{12}^{(n)}} \sum_{j=1}^g\frac{v_j(P_n)}{u_k-q_j}
\\
- \frac{\Omega(P_k)\Omega(P_n)}{16} \sum_{j=1}^g  v_j(P_n) v_j(P_k)
\left(  \frac{1}{2}\sum_{l=1}^{2g+1} \frac{1}{q_j-u_l}
 +   \frac{1}{u_n-q_j} - \sum_{m=1, m\neq j}^g\frac{1}{q_j-q_m}
 \right)
 \\
 -\frac{\Omega(P_k)\Omega(P_n)}{16} \sum_{j=1}^g  v_j(P_n)\sum_{l=1, l\neq j} \frac{v_l(P_k)}{q_j-q_l} + \frac{1}{2}A_{12}^{(k)}A_{12}^{(n)}.
\end{multline*}
Recall that we want to prove that the right hand side of this equality coincides with that of (\ref{beta_diffeq}), that is we need to prove the following
\begin{multline}
\label{temp_dbeta}
\frac{1}{2}\frac{1}{u_k-u_n} \left(   \frac{A_{12}^{(k)}}{A_{12}^{(n)}}\beta_n^2  -  \frac{A_{12}^{(n)}}{A_{12}^{(k)}} \beta_k^2  \right)
=\frac{\Omega(P_n)}{4}A_{12}^{(k)}\frac{\beta_n}{A_{12}^{(n)}} \sum_{j=1}^g\frac{v_j(P_n)}{u_k-q_j}
\\
- \frac{\Omega(P_k)\Omega(P_n)}{16} \sum_{j=1}^g  v_j(P_n) v_j(P_k)
\left(  \frac{1}{2}\sum_{l=1}^{2g+1} \frac{1}{q_j-u_l}
 +   \frac{1}{u_n-q_j} - \sum_{\alpha=1, \alpha\neq j}^g\frac{1}{q_j-q_\alpha}
 \right)
 \\
 -\frac{\Omega(P_k)\Omega(P_n)}{16} \sum_{j=1}^g  v_j(P_n)\sum_{\alpha=1, \alpha\neq j} \frac{v_\alpha(P_k)}{q_j-q_\alpha} + \frac{1}{2}A_{12}^{(k)}A_{12}^{(n)}.
\end{multline}
To shorten expressions which appear next in the calculation, let us adopt the following notation:
\begin{equation}
\label{sigma_n}
\Sigma_n:=\sum_{j=1}^g \frac{1}{\phi(Q_j)(u_n-q_j) \prod_{\alpha=1,\alpha\neq j}^g(q_j-q_\alpha)}
\end{equation}
and
\begin{equation}
\label{sigma_nk}
\Sigma_{nk}:= \sum_{j=1}^g \frac{1}{\phi(Q_j)(u_n-q_j)(u_k-q_j) \prod_{\alpha=1,\alpha\neq j}^g(q_j-q_\alpha)}.
\end{equation}

Note that every term in equality (\ref{temp_dbeta}) contains the product $\Omega(P_n)\Omega(P_k)$ in the numerator. Moreover, we can divide every term by $A_{12}^{(k)}A_{12}^{(n)}$. This gives
\begin{multline}
\label{rhs}
\frac{1}{2}\frac{1}{u_k-u_n} \left(  \left( \frac{\beta_n}{A_{12}^{(n)}}\right)^2  -  \left(\frac{\beta_k}{A_{12}^{(k)}} \right)^2  \right)
= \Sigma_n\Sigma_{nk}
-\frac{\Omega(P_\infty)}{2} \Sigma_{nk}
\\
- \sum_{j=1}^g\frac{1}{\phi^2(Q_j)(u_n-q_j)(u_k-q_j) \prod_{\alpha=1,\alpha\neq j}^g(q_j-q_\alpha)^2} \left(  \frac{1}{2}\sum_{l=1}^{2g+1} \frac{1}{q_j-u_l}
 +   \frac{1}{u_n-q_j} - \sum_{m=1, m\neq j}^g\frac{1}{q_j-q_m}
 \right)
 \\
  - \sum_{j=1}^g \frac{1}{\phi(Q_j)(u_n-q_j) \prod_{\alpha=1,\alpha\neq j}^g(q_j-q_\alpha)}
  \sum_{l=1, l\neq j}^g \frac{1}{\phi(Q_l)(u_k-q_l) \prod_{\alpha=1,\alpha\neq l}^g(q_l-q_\alpha)(q_j-q_l)} + \frac{1}{2}.
\end{multline}
In the left hand side, due to (\ref{boston}), the difference of squares becomes:
\begin{equation}
\label{lhs}
\frac{1}{2}\Sigma_{nk}\left(\Sigma_n+ \Sigma_{k} - \Omega(P_\infty)\right)
\end{equation}
and the terms with $\Omega(P_\infty)$ cancel out. What remains to show now is the following equality:
\begin{multline}
\label{temp_last}
\frac{1}{2}\Sigma_{nk}\left(\Sigma_n- \Sigma_{k}\right) =
 \sum_{j=1}^g\frac{  \frac{1}{2}\sum_{l=1}^{2g+1} \frac{1}{q_j-u_l}
 +   \frac{1}{u_n-q_j} - \sum_{m=1, m\neq j}^g\frac{1}{q_j-q_m}
 }{\phi^2(Q_j)(u_n-q_j)(u_k-q_j) \prod_{\alpha=1,\alpha\neq j}^g(q_j-q_\alpha)^2} - \frac{1}{2}
 \\
+ \sum_{j=1}^g \frac{1}{\phi(Q_j)(u_n-q_j) \prod_{\alpha=1,\alpha\neq j}^g(q_j-q_\alpha)}
   \sum_{l=1, l\neq j}^g \frac{1}{\phi(Q_l)(u_k-q_l) \prod_{\alpha=1,\alpha\neq l}^g(q_l-q_\alpha)(q_j-q_l)}.
\end{multline}
Let us look at the left hand side of this equality. First, note that it can be written as
\begin{multline}
\label{sigmas}
\frac{u_k-u_n}{2}\Sigma_{nk}^2 = \frac{1}{2} \sum_{j=1}^g \frac{u_k-u_n}{\phi^2(Q_j)(u_n-q_j)^2(u_k-q_j)^2 \prod_{\alpha=1,\alpha\neq j}^g(q_j-q_\alpha)^2}
\\
+\frac{1}{2} \sum_{j=1}^g \sum_{l=1,l\neq j}^g\frac{u_k-u_n}{\phi(Q_j)\phi(Q_l)(u_n-q_j)(u_k-q_j) (u_n-q_l)(u_k-q_l) \prod_{\alpha=1,\alpha\neq j}^g(q_j-q_\alpha)  \prod_{\alpha=1,\alpha\neq l}^g(q_l-q_\alpha)}.
\end{multline}
%
The first sum in the right hand side of this equality is contained in the first sum of the right hand side of (\ref{temp_last}), namely if we single out the terms with $l=n$ and $l=k$ in the numerator in (\ref{temp_last}) together with $1/(u_n-q_j)$, we obtain $\frac{1}{2}\left( \frac{1}{u_n-q_j} - \frac{1}{u_k-q_j}  \right)$ which is the same as $\frac{1}{2} \frac{u_k-u_n}{(u_n-q_j)(u_n-q_j)} $. As for the second sum in (\ref{sigmas}), using the following trivial identity
\begin{equation*}
\frac{u_k-u_n}{(u_n-q_j)(u_k-q_j) (u_n-q_l)(u_k-q_l)}= \frac{1}{q_j-q_l} \left( \frac{1}{(u_n-q_j)(u_k-q_l)} - \frac{1}{(u_k-q_j)(u_n-q_l)} \right)
\end{equation*}
we see that it is equal to the sum in the second line of (\ref{temp_last}).
We have thus reduced (\ref{temp_last}) to
\begin{multline}
 \sum_{j=1}^g\frac{ 1 }{\phi^2(Q_j)(u_n-q_j)(u_k-q_j) \prod_{\alpha=1,\alpha\neq j}^g(q_j-q_\alpha)^2} \left(  \frac{1}{2}\sum_{l=1, l\neq n, k}^{2g+1} \frac{1}{q_j-u_l}
  - \sum_{m=1, m\neq j}^g\frac{1}{q_j-q_m}
 \right) = \frac{1}{2}.
 \end{multline}
Rewriting $\phi^{-2}(Q_j)$ as a product using the local parameter $u-q_j$ near $Q_j$ to evaluate according to (\ref{evaluation}) the differential $\phi$ at the points $Q_j$, we see that this equality involves only rational functions. We prove that it holds identically by considering the sum of residues of the following differential form on the Riemann $u$-sphere
%
%
\begin{equation*}
\frac{\prod_{s=1, s\neq n,k}^{2g+1} (u-u_s) du}{\prod_{m=1}^g(u-q_m)^2}.
\end{equation*}
$\Box$

\begin{lemma}
\label{lemma_21}
The $21$-entries of the residue-matrices $A^{(n)}$ \eqref{Ai_hyper} satisfy the following differential equation:
\begin{equation*}
\frac{\partial A_{21}^{n}}{\partial u_k} = \frac{2}{u_k-u_n} \left(  A_{21}^{(k)} A_{11}^{(n)} - A_{21}^{(n)}A_{11}^{(k)}  \right).
\end{equation*}
\end{lemma}

{\it Proof.} As in the elliptic case, see equations \eqref{Schlesinger_21}, this lemma can be seen as a corollary of Lemmas \ref{lemma_diffeq_12} and \ref{lemma_diffeq_11}. Namely, given the form \eqref{Ai_hyper} of $A^{(n)}_{21}$ its derivative can be written as follows
\begin{equation*}
\frac{\partial A_{21}^{n}}{\partial u_k} = -  \frac{1+2\beta_n}{4A_{12}^{(n)}} \frac{\partial \beta_n}{\partial u_k} +  \frac{\beta_n+\beta_n^2}{4\left(A_{12}^{(n)}\right)^2} \frac{\partial A_{12}^{(n)}}{\partial u_k}.
\end{equation*}
From \eqref{Ai_hyper} we se that $\partial_{u_k}\beta_n= -2 \partial_{u_k}A_{11}^{(n)}$. Now using Lemma \ref{lemma_diffeq_11} for $\partial_{u_k}A_{11}^{(n)}$ and (\ref{diffeq_12_hyperelliptic}) for $\partial_{u_k}A_{12}^{(n)}$, we prove the lemma.
$\Box$

\vskip 0.4 cm
We have thus proved that the matrices (\ref{Ai_hyper}) satisfy the part of the Schlesinger system (\ref{Schlesinger_matrix}) with $j\neq k$.
To prove that the remaining differential equations are also satisfied, it suffices to show that the residue of the form $A(u)du$ at infinity is a constant matrix: $A^{(1)} +\dots+A^{(2g+1)} = -A^{(\infty)} = const.$

\begin{lemma}
\label{lemma_sum21}
The sum of residues at $u=u_n$ for  $n=1, \dots, 2g+1$ of the low off-diagonal term of the matrix $A(u)du$ (\ref{Ah}) vanishes:
\begin{equation*}
\sum_{n=1}^{2g+1} A_{21}^{(n)} =0.
\end{equation*}
\end{lemma}

{\it Proof.}
The $21$-entries of the matrices $A^{(n)}$ are defined by (\ref{Ai_hyper}) as follows
\begin{equation*}
A_{21}^{(n)}=- \frac{1}{4} \,\frac{\beta_n+\beta_n^2}{A_{12}^{(n)}}
\end{equation*}
where $A_{12}^{(n)}$ are given by (\ref{def_residues}) and $\beta_n$ by (\ref{betas_hyperelliptic}). Using Lemma \ref{lemma_betasum}, Remark \ref{rmk_sum12} and notation (\ref{sigma_n}), we see that we need to compute the following sum:
\begin{equation*}
-4\sum_{n=1}^{2g+1} A_{21}^{(n)} = \sum_{n=1}^{2g+1} \Sigma_n
+ \sum_{n=1}^{2g+1} A_{12}^{(n)}\left( \Sigma_n  \right)^2   - \frac{\Omega(P_\infty)}{2}.
\end{equation*}

A straightforward computation of the square of the sum $\Sigma_n$ in the second term of the right hand side and changing the order of summation brings us to computing the following two sums (with $i\neq j$):
\begin{equation}
\label{sum1}
\sum_{n=1}^{2g+1} \frac{A_{12}^{(n)}}{(u_n-q_j)(u_n-q_i)} =\frac{\phi(Q_j) \prod_{\alpha=1, \alpha\neq j}^g(q_j-q_\alpha)}{q_i-q_j} + \frac{\phi(Q_i) \prod_{\alpha=1, \alpha\neq i}^g(q_i-q_\alpha)}{q_j-q_i}
\end{equation}
and
\begin{equation}
\label{sum2}
\sum_{n=1}^{2g+1} \frac{A_{12}^{(n)}}{(u_n-q_j)^2} =- \phi(Q_j) \prod_{\alpha=1, \alpha\neq j}^g(q_j-q_\alpha)
\left( \alpha_j + 2\sum_{\alpha=1,\alpha\neq j}^g \frac{1}{q_j-q_\alpha} - \sum_{i=1}^{2g+1} \frac{1}{q_j-u_i}    \right)
\end{equation}
where the $\alpha_j$ in (\ref{sum2}) are the normalization coefficients defined by (\ref{hOmega_vj}). The equalities (\ref{sum1}) and (\ref{sum2}) are obtained by equating to zero the sum of residues of the differentials
\begin{equation*}
\frac{A_{12}(u)\;du}{(u-q_j)(u-q_i)} \qquad\mbox{and}\qquad \frac{A_{12}(u)\,du}{(u-q_j)^2}
\end{equation*}
on the $u$-sphere, respectively. Here $A_{12}(u)$ is the function defined by (\ref{A12_hyperelliptic}) having simple poles at the branch points of the hyperelliptic curve.

Plugging in results (\ref{sum1}) and (\ref{sum2}) and using (\ref{Omega_infty}), we prove the lemma.
$\Box$

Collecting results of Remark \ref{rmk_sum12}, Corollary \ref{cor_sum11} and Lemma \ref{lemma_sum21}, we obtain the matrix $A^{(\infty)}$ (see (\ref{hSchlesinger_intro})).
\begin{corollary}
\label{cor_Ainfty}
The residue at $u=\infty$ of the form $A(u)du$ with the matrix $A$ given by (\ref{Ah}) and (\ref{Ai_hyper}) is a constant matrix:
\begin{equation*}
A^{(\infty)}=-\sum_{n=1}^{2g+1} A^{(n)} = \left(\begin{array}{rr} \frac{1}{4} & 0 \\0 & -\frac{1}{4}\end{array}\right).
\end{equation*}
\end{corollary}
This completes the proof of Theorem \ref{thm_main}.

\begin{remark} {\rm For a generic pair $(\mathcal L, z_0)$, where $\mathcal L$ is a hyperelliptic Riemann surface and $z_0$ is a point in the Jacobian of $\mathcal L$, which satisfies assumptions (i)-(iii) from Section \ref{subsection_inversion}, all zeros of the differential $\Omega$ are simple. In this case,
the $2g-1$ functions $z_i(u_1, \dots, u_{2g-1})$ given by projections of the zeros of the differential $\Omega$ from the hyperelliptic curve (\ref{hyperelliptic}) on the $u$-sphere together with the $2g-1$ functions given by
\begin{equation*}
\frac{1}{4}\sum_{i=1}^{2g-1}\frac{1+2\beta_i}{u_i-z_k}, \qquad k=1, \dots, 2g-1,
\end{equation*}
satisfy the multidimensional Garnier system of the form given by (4.1.9) of \cite{Jap}.
 This is the statement of Corollary $6.2.2$ in \cite{Jap} and is explained as follows. The linear system \eqref{linsys_intro} for $2\times 2$ matrices can be rewritten as a second order ordinary differential equation, see (6.1.2) and (6.1.3) of \cite{Jap}:
\begin{equation}
\label{ode}
\frac{d^2Y}{du^2} + p_1(u; u_1,\dots, u_{2g+1}) \frac{dY}{du} + p_2(u; u_1,\dots, u_{2g+1}) Y=0
\end{equation}
with the coefficients $p_1$ and $p_2$ which, in our case of $A\in sl(2,{\mathbb C})$, are given by
\begin{equation*}
p_1=-\frac{d}{du} {\rm ln} \,A_{12}(u),\qquad
p_2={\rm det} \,A(u) - \frac{d}{du} A_{11} + A_{11} \frac{d}{du} {\rm ln}\, A_{12}(u).
\end{equation*}
The multidimensional Garnier system (4.1.9) from \cite{Jap} is the  condition of isomonodromy deformation of this second order ordinary differential equation. Therefore it is equivalent to the Schlesinger system which expresses the isomonodromy deformation condition for linear matrix system \eqref{linsys_intro}.

Historically the Garnier system appeared as follows. In 1907, R. Fuchs discovered \cite{Fuchs} the Painlev\'e VI equation in the study of isomonodromic deformation condition of a second order linear ordinary differential equation with four essential singularities at $0,1,x, \infty$ and one apparent singularity at $y$. He found out that the monodromy group of the equation is independent of the position  of the singularity $x$, if the position of the apparent singularity $y$ as a function of $x$ satisfies the Painlev\'e's sixth equation.  In 1912, in  paper \cite{Garnier}, Garnier extended this result to the case of more general second order linear ordinary differential equations having $n+3$ essential and $n$ apparent singularities. Thus solutions of the original Garnier system give the positions of apparent singularities as functions of positions of essential singularities such that the monodromy group of the equation is independent of the positions of the essential singularities.

The multidimensional Garnier system studied in \cite{Jap}, is obtained by generalizing the approach of Fuchs and Garnier to the case of a second order ordinary differential equation of the form \eqref{ode} with $n+3$ regular  and $n$ apparent singularities.

For $g=1$ both, the original Garnier and the multidimensional Garnier systems, turn into Painlev\'e's sixth equation, see \cite{Garnier} and \cite{Jap}.

}
\end{remark}

\section{Global behaviour of solutions}
\label{sect_independence}

\subsection{Independence of the choice of homology basis}

In this section we show that, for a fixed elliptic or hyperelliptic surface, the differential $\Omega$ is independent of the choice of a canonical homology basis on the surface.

Let us fix the positive divisor $D=Q_1+\dots + Q_g$ on a hyperelliptic curve of genus $g$.
 Denote by ${\bf a}$ and $\bf b$ the vectors composed of the elements of a canonical homology basis: ${\bf a}=(a_1, \dots, a_g)^t$ and ${\bf b}=(b_1, \dots, b_g)^t$. Any other canonical homology basis is related to $\bf(a,b)$ by a linear transformation with a symplectic matrix:
\begin{equation*}
\left(\begin{array}{c}\tilde{\bf b} \\\tilde{\bf a}\end{array}\right) = \left(\begin{array}{cc}A & B \\C & D\end{array}\right) \left(\begin{array}{c}{\bf b} \\{\bf a}\end{array}\right) \qquad \mbox{with}\qquad
\left(\begin{array}{cc}A & B \\C & D\end{array}\right)\in {\rm Sp}(2g, \mathbb Z).
\end{equation*}
The two corresponding vectors $\omega$ and $\tilde{\omega}$ of holomorphic normalized differentials and the two Riemann matrices $\mathbb B$ and $\tilde{\mathbb B}$  are related by
\begin{equation*}
{\bf \omega}=(C\mathbb B+D)^t\tilde{\bf \omega} \qquad\mbox{and} \qquad
\tilde{\mathbb B}=(A\mathbb B+B)(C\mathbb B+D)^{-1}.
\end{equation*}
%

Since the Abel  map and the Jacobian depend on the choice of a homology basis, the coordinate vectors $c_1$ and $c_2$ (\ref{inversion})  of the Abel image ${\mathcal A}(D)$ of a fixed divisor $D$ also depend on the choice of canonical homology basis. We have:
\begin{equation*}
{\mathcal A}(D)=\int_{P_\infty}^{Q_1} {\bf \omega} + \dots + \int_{P_\infty}^{Q_g} {\bf \omega} = c_1+{\mathbb B}c_2, \qquad \mbox{and} \qquad
\tilde{\mathcal A}(D)= \int_{P_\infty}^{Q_1} { \tilde{\bf \omega}} + \dots + \int_{P_\infty}^{Q_g} \tilde{\bf \omega} = \tilde c_1+\tilde {\mathbb B} \tilde c_2,
\end{equation*}
\begin{lemma}
\label{lemma_ctransformation}
The coordinate vectors $(c_1,c_2)$ and $(\tilde c_1,\tilde c_2)$ are related by
\begin{equation*}
\tilde{c}_1=Ac_1-Bc_2 \qquad \mbox{and} \qquad \tilde{c}_2=Dc_2-Cc_1,
\end{equation*}
where $A,B,C,D$ are $g\times g$ matrices forming the symplectic matrix of transformation between two canonical homology bases $({\bf a}, {\bf b})$ and $(\tilde{\bf a}, \tilde{\bf b})$.
\end{lemma}

{\it Proof.}
Using the transformation rule for the holomorphic normalized differential under the change of basis cycles, for the Abel image $\tilde{\mathcal A}(D)$  of the divisor D, we have
\begin{equation*}
\tilde{\mathcal A}(D)= ((C\mathbb B+D)^{-1})^t{\mathcal A}(D).
\end{equation*}
Therefore, we need to find $\tilde{c}_1$ and $\tilde{c}_2$ such that
\begin{equation}
\label{Abel_transform_1}
 \tilde c_1+\tilde {\mathbb B} \tilde c_2= ((C\mathbb B+D)^{-1})^t(c_1+{\mathbb B}c_2).
\end{equation}
For the  vectors $\tilde c_1$, $\tilde c_2$ from the lemma, we get
\begin{equation*}
 \tilde c_1+\tilde {\mathbb B} \tilde c_2 = (A-\tilde{\mathbb B}C)c_1+(\tilde{\mathbb B}D-B)c_2.
\end{equation*}
Now let us look at the terms $\tilde{\mathbb B}C$ and $\tilde{\mathbb B}D$. Since the Riemann matrix is symmetric, we can rewrite the formula relating $\tilde{\mathbb B}$ and $\mathbb B$ as follows: $\tilde{\mathbb B}=((C\mathbb B+D)^{-1})^t(A\mathbb B+B)^t$. Because the matrices $A, B, C, D$ form a symplectic matrix, they satisfy $D^tA-B^tC=\mathbb I$ and $C^tA=A^tC.$ Using this, we get $\tilde{\mathbb B}C=A-((C\mathbb B+D)^{-1})^t.$
%
%
Similarly, using one more relation for blocks of a symplectic matrix, namely $D^tB=B^tD$, we obtain  $\tilde{\mathbb B}D=((C\mathbb B+D)^{-1})^t\mathbb B+B.$ This proves that (\ref{Abel_transform_1}) holds for the vectors $\tilde c_1$, $\tilde c_2$ from the statement of the lemma.
$\Box$

\begin{proposition}
The differential $\Omega$ is independent of the choice of canonical homology basis on a hyperelliptic Riemann surface.
\end{proposition}
{\it Proof.}
We defined $\Omega$ as an Abelian differential of the third kind with simple poles at the divisor $D+D^*$ with residues $+1$ at $Q_j$ and $-1$ at $Q^*_j$ normalized to have the vector of $a$-periods equal to $-4\pi{\rm i}c_2$. Here $(c_1, c_2)$ are the coordinates of the Abel image of $D$ with respect to the lattice of the Jacobian.
 As we change the canonical homology basis, the $a$-cycles will transform according to $\tilde{\bf a}=C{\bf b}+D{\bf a}.$ Thus the vector of $a$-periods of $\Omega$ with respect to the new homology basis is given by $4\pi{\rm i}Cc_1-4\pi{\rm i}Dc_2$. This vector, due to Lemma \ref{lemma_ctransformation},  coincides with $-4\pi{\rm i}\tilde{c}_2$, where $(\tilde{c}_1, \tilde{c}_2)$ are the coordinates of the Abel image of $D$ with respect to the lattice of the Jacobian after the transformation of the basis of cycles.
 Therefore, the differential $\Omega$ satisfies the correct normalization condition with respect to the new basis of cycles $\tilde{\bf a}, \tilde{\bf b}$.
$\Box$

\subsection{Analytical continuation of solutions}

Now let us consider the situation described in Sections \ref{sect_EllipticOkamoto} and \ref{sect_ESchlesinger}, that is the family of two-fold elliptic coverings. The results presented there, namely the description of the general solution $y(x)$ of the Painlev\'e VI with parameters (\ref{constants}) and of the solution of the corresponding Schlesinger system, were obtained for small variations of the moving branch point $x$.

Let us now consider an analytical continuation of the general solution $y(x)$ as the moving branch point $x$ makes a full tour around another branch point and comes back to its initial position. Recall that we define the point $Q_0$ as the Jacobi inversion of the point $z_0$ having constant coordinates with respect to the lattice $\Lambda$ generated by the periods $1$ and $\mu$ of the elliptic curve: $z_0=c_1+c_2\mu$.
As the point $x$ returns to its original position after going around another branch point, the $a$- and $b$-cycles of the surface transform in a non-trivial way. Thus the lattice $\Lambda$ also transforms. Since the coordinates $c_1$ and $c_2$ are kept unchanged, it follows from  Lemma \ref{lemma_ctransformation} that as the branch point $x$ comes back after going around another branch point, the point $Q_0$ does not return to its initial position on the surface. This shows that the Picard solution $y_0(x)$ \eqref{Picard} is a multivalued meromorphic function which has monodromies at points $0,1$ and $\infty$.

Our interpretation of the solution $y(x)$ of the Painlev\'e's sixth equation with parameters \eqref{constants} as the position of zeros of $\Omega$ shows that $y(x)$ can be analytically continued everywhere in the Riemann sphere except the points $x=0,1,\infty.$ The above discussion implies that $y(x)$ has monodromy at these points and only at these points. It is also clear that there are no essential singularities.  Thus $y(x)$ is a multivalued meromorphic function on $\mathbb C\setminus\{0,1\}$. In other words, $y(x)$ enjoys the Painlev\'e property.

The statement about the Painlev\'e property is proved in \cite{HinkkanenLaine} and \cite{OkamotoTakano} for solutions of all Painlev\'e VI equations.

As for the results of Section \ref{sect_HSchlesinger}, the generalization to hyperelliptic curves, in order to understand the global behaviour of our solution, we need to prove that assumptions (i)-(iii) for the divisor $D$, see Subsection \ref{subsection_inversion}, are preserved under our deformation of the pair: hyperelliptic curve and a point $z_0$ in its Jacobian. This remains an open question.

\section{Tau function of the Schlesinger system}
\label{sect_tau}

With any monodromy preserving deformation of a system of linear ordinary differential equations one associates, as described in \cite{JMU}, the so-called {\it isomonodromic tau-function} $\tau$, a function of deformation parameters playing an important role in the theory of isomonodromic deformations. It was proved in \cite{Bolibruch, Malgrange} that the tau-function is holomorphic everywhere in the space of deformation parameters outside of the hyperplanes where the values of  two deformation parameters coincide.

In our case, the Schlesinger system (\ref{hSchlesinger_intro}) describes isomonodromic deformations of linear system (\ref{linsys_intro}), the deformation parameters being the branch points $\{u_j\}$ of the hyperelliptic covering.
The corresponding tau-function is thus holomorphic on the universal covering of the space ${\mathbb C}^n\setminus \{(u_1,\dots,u_{2g+1})\mid u_k=u_l {\mbox{ with }} k\neq l\}.$  The set of zeros of the function $\tau$ in this space is called  {\it the Malgrange divisor}; a solution to the Schlesinger system always exists outside of the Malgrange divisor \cite{BolibruchOrders}.

The tau-function of the Schlesinger system is defined, up to a constant factor, by the following differential equations:
\begin{equation}
\label{tau-def}
\frac{\partial\;{\rm ln} \tau}{\partial u_j} = \frac{1}{2}\;\underset{u=u_j}{\rm res} {\rm tr}\; A^2(u)\,du.
\end{equation}

The tau-functions of the solutions of the Schlesinger systems constructed in this paper coincide with those from \cite{KiKo}. In this section, we prove this statement for the case of the Schlesinger system (\ref{Schlesinger}) corresponding to elliptic coverings with branch points at $\{0,1,x,\infty\}$. In other words, the following theorem holds.

\begin{theorem}
\label{thm_tau}
Consider  solution (\ref{betas})-(\ref{Ai}) from Theorem \ref{thm_Schlesinger} of the Schlesinger system associated with the family of elliptic curves ramified over the set $\{0,1,x,\infty\}$. Its tau-function
is given by
\begin{equation}
\label{tau_elliptic}
\tau(x) = C \frac{\theta[c_2+\frac{1}{2}, c_1+\frac{1}{2}](0)}{\left( x(x-1)\right)^{1/8}\sqrt{I_0}},
\end{equation}
where $C\in {\mathbb C}$ is a constant; $c_1, c_2\in \mathbb C$ are the constants from (\ref{z0}) corresponding to the solution of the Schlesinger system such that $[c_1, c_2]$ is not a half-integer characteristic, and $I_0$ is the $a$-period of the holomorphic differential as in (\ref{I0}).
\end{theorem}

We thus see that, in genus one, the tau-function is holomorphic everywhere outside of the set $x\in\{0,1, \infty\}$ and does not have zeros, therefore the Malgrange divisor is empty.

The rest of the section is devoted to a proof of Theorem \ref{thm_tau} . We start by proving two lemmas.

\begin{lemma}
\label{lemma_Omega}
Let $\mu$ be the period of the elliptic curve (\ref{ell_curve}) and $I_0$ be the elliptic integral defined by (\ref{I0}). The following relation holds.
\begin{equation*}
\frac{\omega^2(P_x)}{16}\left(   \frac{\Omega(P_\infty)}{\omega(P_\infty)}   \right)^2  = \frac{y_0-x}{4x(x-1)} + \frac{d}{dx}\left(  \pi{\rm i} c_2^2 \mu - \frac{1}{2}{\rm ln} \;I_0 + {\rm ln} \;\theta_1\left( \int_{Q_0}^{P_\infty} \omega \right)   \right).
\end{equation*}
\end{lemma}

{\it Proof}.
To obtain this relation we  rewrite the differential $\Omega$ in terms of theta-functions corresponding to the elliptic curve. The fundamental bidifferential $W$ from Section \ref{sect_Hurwitz} can be written as follows:
\begin{equation}
\label{Wtheta}
W(P,Q) = d_P d_Q {\rm ln}\; \theta_1\left( \int_Q^P \omega \right) =
\left[ \left( \frac{\theta_1'\left( \int_Q^P \omega \right)}{\theta_1\left( \int_Q^P \omega \right)}   \right)^2 - \frac{\theta_1''\left( \int_Q^P \omega \right)}{\theta_1\left( \int_Q^P \omega \right)}\right] \omega(P ) \omega(Q).
\end{equation}
The differential of the third kind with poles at $Q_0$ and $Q_0^*$ normalized by vanishing of the $a$-period  is the integral of $W(P,Q)$ over a path from $Q_0^*$ to  $Q_0$ lying in a fundamental polygon of the curve, that is not crossing the basis $a$- and $b$-cycles. From the definition (\ref{Omega}) of $\Omega$ and taking into account that $\theta_1$ is odd and $\int_{Q_0}^{P_\infty}\omega = -\int_{Q_0^*}^{P_\infty}\omega$, we have:

\begin{equation*}
\Omega(P_\infty ) =\left( 2  \frac{\theta_1'\left( \int_{Q_0}^{P_\infty} \omega \right)}{\theta_1\left( \int_{Q_0}^{P_\infty} \omega \right)} -4\pi{\rm i}c_2 \right)\omega(P_\infty).
\end{equation*}
Thus the left hand side  of the relation from lemma becomes
\begin{equation*}
\frac{\omega^2(P_x)}{16}\left(   \frac{\Omega(P_\infty)}{\omega(P_\infty)}   \right)^2 = \frac{\omega^2(P_x)}{4}
\left[ \left(\frac{\theta_1'\left( \int_{Q_0}^{P_\infty} \omega \right)}{\theta_1\left( \int_{Q_0}^{P_\infty} \omega \right)}\right)^2
- 4\pi{\rm i}c_2 \frac{\theta_1'\left( \int_{Q_0}^{P_\infty} \omega \right)}{\theta_1\left( \int_{Q_0}^{P_\infty} \omega \right)}    \right] -{\omega^2(P_x)}\pi^2c_2^2.
\end{equation*}
Adding and subtracting the term ${\theta_1''\left( \int_{Q_0}^{P_\infty} \omega \right)}/{\theta_1\left( \int_{Q_0}^{P_\infty} \omega \right)}$ we single out the presence of the bidifferential $W$ in the form (\ref{Wtheta}) evaluated at the points $P_\infty$ and $Q_0$. On the other hand, the remaining terms with theta-functions make up the derivative $d_x{\rm ln} \;\theta_1\left( \int_{Q_0}^{P_\infty} \omega \right)$ due to definition (\ref{Q0}) of the point $Q_0$ on the elliptic curve, the heat equation for the theta function
\begin{equation*}
\frac{\partial^2\theta_1(z)}{\partial z^2} = 4\pi{\rm i} \frac{\partial \theta_1(z)}{\partial \mu},
\end{equation*}
and the Rauch formula (\ref{Rauchx}) for $d_x\mu$ in genus one. In other words, we have
\begin{equation}
\label{almost}
\frac{\omega^2(P_x)}{16}\left(   \frac{\Omega(P_\infty)}{\omega(P_\infty)}   \right)^2 =
\pi {\rm i}c_2^2\frac{d}{dx}\mu
+\frac{\omega^2(P_x)}{4\omega(Q_0)\omega(P_\infty)}W(Q_0,P_\infty)
+ \frac{d}{dx}{\rm ln}\;\theta_1\left( \int_{Q_0}^{P_\infty} \omega\right).
\end{equation}
Now it remains to integrate the second term in the last line. Note that for the bidifferential $W$ with one argument evaluated at $P_\infty$ we can write
\begin{equation*}
W(P,P_\infty) = \left( u(P )-\frac{I^\infty}{I_0}\right)\frac{\omega(P )}{\omega(P_\infty)},
\end{equation*}
where $u(P )$ is the $u$-coordinate of the point $P$ of the elliptic curve (\ref{ell_curve}) and $I^\infty$ is the elliptic integral of the form $\oint_a u(P )\phi(P ) $. Using this, we can transform $W(Q_0, P_\infty)$ into $W(P_x, P_\infty)$  as follows:
\begin{equation*}
\frac{\omega(P_\infty)}{\omega(Q_0)} W(Q_0,P_\infty) = W(P_x, P_\infty) \frac{\omega(P_\infty)}{\omega(P_x)} + y_0-x.
\end{equation*}
Now using a similar expression (\ref{WPPx}) for the $W$ with one argument being evaluated at $P_x$  with $P=P_\infty$ we have
\begin{equation*}
W(P_\infty,P_x) = -\frac{I^x}{I_0} \frac{\omega(P_\infty)}{\omega(P_x)}.
\end{equation*}
Combining all this, we get the following for the second term in (\ref{almost}) which is the last term we need to integrate to prove the lemma:
\begin{multline*}
\frac{\omega^2(P_x)}{4\omega(Q_0)\omega(P_\infty)}W(Q_0,P_\infty) = \frac{\omega^2(P_x)}{4\omega^2(P_\infty)} \left(  W(P_x, P_\infty) \frac{\omega(P_\infty)}{\omega(P_x)} + y_0-x  \right) = -\frac{I^x}{4I_0} + \frac{y_0-x}{4x(x-1)},
\end{multline*}
where we used (\ref{omega_atPx}) for $\omega(P_x)$ and that $\omega(P_\infty) = -{2}/{I_0}$ as can be seen by computing the evaluation defined by (\ref{evaluation}) of the holomorphic normalized differential $\omega$ (\ref{holdiff}) at the point $P_\infty$ with respect to the local parameter given by  ${u}^{-1/2}$.
%
 From definitions (\ref{I0}), (\ref{Ix}) of the integrals $I_0$ and $I^x$, we derive that $I^x$ is the $x$-derivative of $2I_0$, which finishes the proof.

$\Box$

 \begin{lemma}
 \label{lemma_theta}
 Let as before $\mu$ be the period of the elliptic curve. The following relation holds for the theta-functions corresponding to the curve:
 \begin{equation*}
 \theta_1\left( \int_{Q_0}^{P_\infty} \omega \right) {\rm e}^{\pi{\rm i} c_2^2 \mu} = \theta\left[c_2+\frac{1}{2}, c_1+\frac{1}{2}\right](0)\;{\rm e}^{-2\pi{\rm i}c_1c_2-\pi\rm{i}c_2},
 \end{equation*}
 where $\left[c_2+\frac{1}{2}, c_1+\frac{1}{2}\right]$ is a characteristic of the theta-function.
 \end{lemma}

{\it Proof}. By definition of the point $Q_0$ on the elliptic curve, see (\ref{Q0}), the argument of the $\theta_1$ in the left hand side is equal to $-c_1-c_2\mu$. Note that by definition $\theta_1(z) = -\theta[1/2,1/2](z)$ and the characteristic of the theta-function can be transformed into the argument and vice versa. Thus the rest of the proof is a simple calculation using the following relation immediately following from the properties of theta-functions:
 \begin{equation*}
 \theta[p,q](\alpha\mu+\beta) = \theta[\alpha+p, \beta+q](0)\; {\rm e}^{-\pi{\rm i}\alpha^2\mu - 2\pi{\rm i} \alpha\beta - 2\pi{\rm i} \alpha q}.
 \end{equation*}
 $\Box$

{\it Proof of Theorem \ref{thm_tau}.}
For the Schlesinger system (\ref{Schlesinger}) corresponding to an elliptic covering with branch points at $\{0,1,x,\infty\}$,
 the right-hand side of (\ref{tau-def}) can be rewritten as follows

\begin{equation*}
\frac{1}{2}\underset{u=x}{\rm res}\, {\rm tr} A^2(u)du =
{\rm tr} \left(  \frac{A^{(1)}A^{(3)}}{x}  + \frac{A^{(2)}A^{(3)}}{x-1}   \right).
\end{equation*}
%
%
%
Let us compute the traces and write the result using the matrices $A^{(i)}$ of our solution (\ref{Ai}) of the Schlesinger system (\ref{Schlesinger}). Putting $u_1=0$ and $u_2=1$, and after some simple rearrangement of terms, we get
%
%

%
\begin{equation*}
\frac{d\;{\rm ln} \tau}{dx}  = \sum_{k=1}^2 \frac{1}{x-u_k} \left\{ \frac{1}{8} + \frac{1}{4} (\beta_k+\beta_3) - \frac{1}{4} A_{12}^{(k)}A_{12}^{(3)}\left(  \frac{\beta_3}{A_{12}^{(3)}}-\frac{\beta_k}{A_{12}^{(k)}}\right)^2 - \frac{1}{4} \frac{A_{12}^{(k)}} {A_{12}^{(3)}}\beta_3 -\frac{1}{4} \frac{A_{12}^{(3)}} {A_{12}^{(k)}}\beta_k \right\}.
\end{equation*}
Now expressing $\beta_k$ in terms of $A_{12}^{(k)}$ using their definition (\ref{betas}), (\ref{A12_residues})  similarly to (\ref{beta1_A12}) and using formulas (\ref{Omega_0x}), (\ref{Omega_1x}), we have

\begin{multline*}
\frac{d\;{\rm ln} \tau}{dx}  = \frac{1}{8}\left( \frac{1}{x}+\frac{1}{x-1}\right)
+ \frac{ x( \beta_1 +\beta_3 +\beta_2) -\beta_1-\beta_3+\beta_3x  }{4x(x-1)}
-\frac{1}{4} \left( A_{12}^{(1)}A_{12}^{(3)} \frac{x(y_0-1)}{y_0(y_0-x)} \right. \\
\left.- A_{12}^{(1)}A_{12}^{(3)}\frac{x-1}{x-y_0}
+ \frac{1}{y_0}A_{12}^{(1)}A_{12}^{(3)} - A_{12}^{(2)}A_{12}^{(3)} \frac{y_0(x-1)}{(y_0-1)(y_0-x)} - A_{12}^{(2)}A_{12}^{(3)}\frac{x}{x-y_0} - A_{12}^{(2)}A_{12}^{(3)}\frac{1}{y_0-1}\right).
\end{multline*}
Here all the terms with the product  $A_{12}^{(k)}A_{12}^{(3)}$, $k=1,2$ cancel out. In the first line we use Lemma \ref{lemma_betasum} which shows that the sum of the three $\beta_k$ is equal to $-1$ thus reducing the expression to:
\begin{equation}
\label{temp_dlntau}
\frac{d\;{\rm ln} \tau}{dx}  = \frac{1}{8}\left( \frac{1}{x}-\frac{1}{x-1}\right)
 - \frac{\beta_1}{4x(x-1)} + \frac{\beta_3}{4x}.
\end{equation}
According to their definition (\ref{betas}), the coefficients $\beta_1, \beta_3$ contain the quantities $\frac{\Omega(P_0)}{\phi(P_0)}$ and $\frac{\Omega(P_x)}{\phi(P_x)}$, respectively,
where following the notation from Section \ref{sect_HSchlesinger}, $\phi(P )=I_0\omega(P )$. We want to express these quantities in terms of $\frac{\Omega(P_\infty)}{\phi(P_\infty)}$. Note that due to (\ref{Omega_omega}),
%
%
%
\begin{equation*}
\frac{\Omega(P_0)}{\phi(P_0)} = -\frac{1}{y_0\phi(Q_0)} +\frac{\Omega(P_\infty)}{\phi(P_\infty)}  \qquad \mbox{and} \qquad
\frac{\Omega(P_x)}{\phi(P_x)} = \frac{1}{(x-y_0)\phi(Q_0)} + \frac{\Omega(P_\infty)}{\phi(P_\infty)}.
\end{equation*}
Substituting this into (\ref{temp_dlntau}), we have
\begin{equation*}
\frac{d\;{\rm ln} \tau}{dx}  = \frac{1}{8}\left( \frac{1}{x}-\frac{1}{x-1}\right)
+ \frac{1}{4x(x-1)}\left( \frac{\Omega(P_\infty)}{\phi(P_\infty)} \right)^2 - \frac{y_0-1}{4x(x-1)}.
\end{equation*}
Now rewriting ${1}/{4x(x-1)}$ as $\phi^2(P_x)/16$ and using the results of Lemmas \ref{lemma_Omega} and \ref{lemma_theta}, we complete the proof of the theorem. Note that the exponential factor independent of the variable $x$ is unessential because the tau-function is only defined up to multiplication by a constant.

$\Box$

\vskip 2cm

\section{Conclusion: back to Poncelet. The Billiard Ordered Games}

\label{sect_backToPoncelet}

\subsection{Mechanical interpretation of the Poncelet theorem}

When we discussed the elliptic case, we observed, following
Hitchin \cite{Hitchin}, that the integrality condition
(\ref{kJacobian}) has been connected to the Poncelet polygons and
the Cayley condition \cite{Cayley}.

A natural question is whether there is an analogue of such
geometric background in the hyperelliptic cases as well. By
analogy with the elliptic case, we should restrict ourselves to
the assumption that in the  conditions (\ref{preinversion}),
 (\ref{inversion}), the vectors $c_1, c_2$ have all rational
components. Denote by  $\Lambda$ the lattice of the Jacobian of the hyperelliptic curve of
genus $g$. Relation (\ref{inversion}) for rational constant vectors
$c_1, c_2$ can be equivalently reformulated as

\begin{equation}
\label{inversionrat} \int_{P_\infty}^{Q_1} {\bf \omega} + \dots +
\int_{P_\infty}^{Q_g} {\bf \omega} \in \frac{1}{n}\Lambda.
\end{equation}

Our goal now is to present a higher-dimensional generalization of
the Poncelet polygons, for which the Cayley-type condition of
periodicity is equivalent to relation
(\ref{inversionrat}). The higher-dimensional generalizations of
the Poncelet polygons with the required property appear to be the
periodic trajectories of a special class of the billiard ordered
games. The billiard ordered games were introduced recently in
\cite {DR2006} as a generalization of billiards in a conic.

Recall that the Poncelet polygons, after an affine transformation,
can be seen as periodic trajectories of billiards in  a conic.
We will review this mechanical side of the  Poncelet polygons in
more detail now.

\emph{The elliptical billiard} \cite{KT1991}, \cite{Tab2005book}
is a dynamical system where a material point of the unit mass is
moving under the inertia, or in other words, with a constant
velocity inside an ellipse and obeying the reflection law at the
boundary, i.e.\ having congruent impact and reflection angles with
the tangent line to the ellipse at any bouncing point.

Reflections are assumed to be absolutely elastic and
friction is neglected.

It is well known that all segments of a given elliptical billiard
trajectory are tangent to the same conic, the caustic, which is
confocal with the boundary \cite{CCS1993}.

If a trajectory becomes closed after $n$ reflections, then it
presents a Poncelet polygon. The Poncelet Theorem implies that in this case any
trajectory of the billiard system which shares the same caustic
is also periodic with the period $n$. A more general statement,
the Full Poncelet Theorem also has a mechanical meaning.
We will formulate it in the so-called dual version, because
this is the form we need for higher-dimensional generalizations:
suppose we are given $n+1$ conics  in $\mathbb R^2$ from a confocal
pencil, a caustic and $n$ boundary conics. If there exists an $n$-periodic
billiard trajectory with the given caustic and with vertices of
the billiard reflections at each of $n$ boundaries, then there exists such an $n$-periodic billiard
trajectory having a vertex at an arbitrary point of any of the $n$
given boundaries.

The statement dual to the Full Poncelet Theorem can be generalized
to the $d$-di\-men\-sio\-nal space \cite{CCS1993} (see also
\cite{Previato2002}). Darboux  investigated light-rays in the
three-dimensional ($d=3$) case and announced the corresponding
Full Poncelet theorem in \cite{Darboux1870} in 1870.  A systematic
exposition of this higher-dimensional theory is presented
in the book \cite{DR2011}.

\

\subsection{Billiards within Quadrics: $d$-dimensional Case. Billiard ordered games}

\

We are considering billiards in domains of $\mathbb R^d$ bounded by parts of
several confocal quadrics from the family:
\begin{equation}\label{eq:confocal1}
\mathcal{Q}_{\lambda}\ :\ \frac{x_1^2}{a_1-\lambda}+\dots
+\frac{x_d^2}{a_d-\lambda}=1, \quad\lambda\in\mathbb{R},
\end{equation}
with $a_1>a_2>\dots > a_d> 0$ being constants. The Jacobi coordinates in $\mathbb R^d$ are defined as follows. Any point of the
space is contained in exactly  $d$ quadrics from the family
(\ref{eq:confocal1}). The parameters $\lambda_1,\dots,$
$\lambda_d$ of these quadrics are the Jacobi coordinates of the point in $\mathbb R^d$
\cite{Jac}. We always assume
$\lambda_1>\lambda_2>\dots>\lambda_d$.

The  Chasles theorem  states that any line in  $\mathbb
R^d$ is tangent to exactly $d-1$ quadrics from a given confocal
family. The billiard flow has $d-1$ quadrics from the family
(\ref{eq:confocal1}) as caustics. We will denote by $\alpha_1,
\dots, \alpha_{d-1}$ the values of the parameter $\lambda$ which
correspond to the caustics.

Billiards in domains bounded by several confocal quadrics in
arbitrary dimension have been studied and  their periodic
trajectories have been described, see \cite{DR2006},
\cite{DR2011}. The dynamics of such a billiard system can be
studied on the Jacobian of the hyperelliptic curve:
$$
\nu^2=\mathcal P_{2d-1}(\lambda)=(\lambda-a_1)\cdots(\lambda-a_d)(\lambda-\alpha_1)\cdots
(\lambda-\alpha_{d-1}).
$$
The generalized Cayley-type conditions have been derived there
from the study of points of finite order on the Jacobian.

\smallskip

The main goal of this section is to review the Cayley-type conditions
for the generalized dual Full Poncelet Theorem which correspond to the
so-called {\it billiard ordered games}. In general, a billiard
ordered game in $d$-dimensional space can be defined for an
arbitrary number $k$ of  quadrics. However, only a special class
will be of interest for us here: the periodic billiard ordered
games in $d$-dimensional space for $k=d-1$ quadrics will
provide the geometric situation related to  the class of solutions
of the Schlesinger system constructed above from hyperelliptic
curves of genus $g=d-1$ exactly in the same way as the Poncelet
polygons corresponded to the solutions of the Painlev\'e VI
equation in the Hitchin construction for $d=2, \;\;g=1$.

\smallskip

Following \cite{DR2006}, we are going to define a notion of
bounces ``from outside'' and ``from inside''. Let us consider a
quadric $\mathcal Q_{\lambda}$ from the confocal family
(\ref{eq:confocal1}) such that $\lambda \in (a_{s+1}, a_s)$ for
some $s\in \{1,\dots, d\}$, where $a_{d+1}=-\infty$. Observe that
along a billiard ray which reflects at $\mathcal Q_{\lambda}$, the
Jacobi coordinate $\lambda_s$ has a local extremum at the point of
reflection.


A ray reflects \emph{from outside}
at the quadric $\mathcal Q_{\lambda}$ if the reflection point is a
local maximum of the Jacobi coordinate $\lambda_s$, and it
reflects \emph{from inside} if the reflection point is a local
minimum of the coordinate $\lambda_s$.

\smallskip

Assume now that a $(d-1)$-tuple of confocal quadrics $\mathcal
Q_{\beta_1},\dots, \mathcal Q_{\beta_{d-1}}$ from the confocal pencil
(\ref{eq:confocal1}) is given, and $(i_1,\dots,i_{d-1})\in\{-1,1\}^{d-1}$.

By definition, \emph{the billiard ordered game} joined to quadrics
$\mathcal{Q}_{\beta_1}$, \dots, $\mathcal{Q}_{\beta_{d-1}}$, with the
\emph{signature} $(i_1,\dots,i_{d-1})$ is the billiard system with
trajectories having bounces at $\mathcal
Q_{\beta_1},\dots,\mathcal Q_{\beta_{d-1}}$, such that
$$
\displaylines{
    \text{the reflection at}\; \mathcal Q_{\beta_s}\;
    \text{is from inside if}\; i_s=+1;
\cr
    \text{the reflection at}\; \mathcal Q_{\beta_s}\;
    \text{is from outside if}\; i_s=-1.
}
$$

Note that any trajectory of a billiard ordered game has $d-1$
caustics from the same family (\ref{eq:confocal1}).

\smallskip

Suppose now that $\mathcal Q_{\beta_1}$, \dots, $\mathcal
Q_{\beta_{d-1}}$ are ellipsoids, meaning that $\beta_i<a_d, \,
i=1,\dots, d-1$.  Consider a billiard ordered game with the
signature $(i_1,\dots,i_{d-1})$. In order for the trajectories of
such a game to stay bounded, in addition to the assumption that
$\mathcal Q_{\beta_i}$ are ellipsoids, the following condition
should be satisfied:
$$
i_s=-1\ \Rightarrow\ i_{s+1}=i_{s-1}=1 \; \text{and}\;
\beta_{s+1}<\beta_s,\ \beta_{s-1}<\beta_s.
$$
We identify indices 0 and $d$ with $d-1$ and 1, respectively. From
\cite{DR2006}, (see also \cite{DR2011}) we get

\begin{proposition}\label{uslov.igra}
Given a billiard ordered game joined to $d-1$ ellipsoids
$\mathcal{Q}_{\beta_1}$, \dots, $\mathcal Q_{\beta_{d-1}}$ with a
signature $(i_1,\dots,i_{d-1})$. Fix the caustics
$\mathcal{Q}_{\alpha_1}$, \dots, $\mathcal Q_{\alpha_{d-1}}$. Then
a trajectory with these caustics is periodic if and only if
$$
\sum_{s=1}^{d-1} \mathcal A(P_{\beta_s})\in \frac{1}{n}\Lambda,
$$
where
$\Lambda$ is the lattice of the
Jacobian  of the curve $\Gamma\, :\, \nu^2=\mathcal P_{2d-1}(\lambda)$,
the point $P_{\beta_s}$ on the curve $\Gamma$ corresponds to the pair
$(\mathcal{Q}_{\beta_s}, i_s)$ as follows
$P_{\beta_s}=(\beta_s, i_s\sqrt{\mathcal P_{2d-1}(\beta_s)})$,
 and $\mathcal A$ is the Abel map based at the point at infinity.
\end{proposition}
\smallskip

The equivalence between the condition in Proposition
\ref{uslov.igra} and relation (\ref{inversionrat}) is obvious.

\smallskip

As an example, we can present the periodicity condition in the
case of the billiard motion for two ellipsoids in the
three-dimensional case.

\begin{example}
The condition that there exists a closed billiard trajectory for
two ellipsoids $\mathcal Q_{\beta_1}$ and $\mathcal Q_{\beta_2}$,
which bounces exactly $m$ times of each of them, with caustics
$\mathcal Q_{\alpha_1}, \mathcal Q_{\alpha_{2}}$, is:
$$
\rank\left(\begin{array}{llll} f_1'(P_{\beta_2}) & f_2'(P_{\beta_2})
& \dots &
f_{m-2}'(P_{\beta_2}) \\
f_1''(P_{\beta_2}) & f_2''(P_{\beta_2}) & \dots &
f_{m-2}''(P_{\beta_2}) \\
 & & \dots\\
 & & \dots\\
f_1^{(m-1)}(P_{\beta_2}) & f_2^{(m-1)}(P_{\beta_2}) & \dots &
f_{m-2}^{(m-1)}(P_{\beta_2})
\end{array}\right)<m-2.
$$
Here
$$
f_j=
\frac{\nu-B_0-B_1(\lambda-\beta_1)-\dots-B_{j+1}(\lambda-\beta_1)^{j+1}}{\lambda^{j+1}},
\quad 1\le j\le m-2,
$$
and $\nu=B_0+B_1(\lambda-\beta_1)+\dots$ is the Taylor expansion of $\nu(\lambda)$ around
the point symmetric to $P_{\beta_1}$ with respect to the
hyperelliptic involution of the genus 2 curve $\Gamma$. The
notation is the same as in Proposition
\ref{uslov.igra}. This condition corresponds to solutions of the
Schlesinger system with $g=2$ with rational constants
$$
c_{1}=(C_{11}/m, C_{12}/m),\quad c_2=(C_{21}/m, C_{22}/m),
$$
where $C_{ij}$ are integer constants for $i,j\in\{1,2\}$.
\end{example}

Thus there is a full parallelism between Hitchin's observation
relating the Poncelet polygons with algebraic solutions of the
Painlev\'e VI equation and the relationship between the periodic trajectories of the
billiard ordered games in $\mathbb R^{g+1}$  for $g$
ellipsoids and the class of solutions of the Schlesinger systems
for rational constants $c_1, c_2$.

\vskip 1cm

{\bf Acknowledgements.} The authors thank D. Korotkin for  useful
discussions. The research has been partially supported by the NSF grant 1444147. The research of the first author has been partially
supported by the grant 174020 ``Geometry and topology of manifolds,
classical mechanics, and integrable dynamical systems" of the
Ministry of Education and Sciences of Serbia and by the University
of Texas at Dallas. The second author gratefully acknowledges
support from the Natural Sciences and Engineering Research Council of Canada, Fonds de recherche du Qu\'ebec Nature et Technologies (grant in the program ``\'Etablissement de nouveaux chercheurs universitaires") and the University of Sherbrooke.

\end{document}